\documentclass[12pt,letterpaper,reqno]{amsart}
\usepackage{amsfonts}
\usepackage{color}
\usepackage{epsfig}
\usepackage{fancyhdr}
\usepackage[latin1]{inputenc}
\usepackage{amsmath}
 \usepackage{mathrsfs}
 \usepackage{enumerate}
\usepackage{natbib}
\usepackage{graphicx}
\usepackage[colorlinks=true,linkcolor=red,citecolor=blue]{hyperref}
\usepackage{amssymb}
\usepackage{dsfont}
\usepackage{fancyhdr}
\usepackage{multirow}

\renewcommand{\P }{\mathbb{P}} 
\newcommand{\E }{\mathbb{E}}

\newcommand{\R }{\mathbb{R}}

\newcommand{\cM }{\mathcal{M}} 
\newcommand{\K}{\mathbb{K}} 
\newcommand{\id}{\mathds{1}} 
\newcommand{\Om}{\mathcal{X} \times \mathcal{Y}} 
\newcommand{\Dp}{D_{\varphi}} 
\newcommand{\Ddp}{\mathcal{D}_{\varphi}} 
\newcommand{\cqfd}{\hfill$\Box$}
\newcommand{\bo}{\boldsymbol{0}} 
\newcommand{\der}{\mbox{d}}

\newtheorem{theorem}{Theorem}[section]

\newtheorem{proposition}[theorem]{Proposition}

\newtheorem{example}[theorem]{Example}
\newtheorem{remark}[theorem]{Remark}

\begin{document}

\title[Semiparametric estimation of mutual information 
and related criteria]{Semiparametric estimation of mutual information 
and related criteria : optimal test of independence}

\author{Amor Keziou$^1$}
\address{$^1$Laboratoire de Math\'ematiques de Reims EA 4535 and 
ARC-Mathématiques CNRS 3399, Université de Reims Champagne-Ardenne, France}
\email{amor.keziou@univ-reims.fr}

\author{Philippe Regnault$^2$}
\address{$^2$Laboratoire de Math\'ematiques de Reims EA 4535 and 
ARC-Mathématiques CNRS 3399, Université de Reims Champagne-Ardenne, France}
\email{philippe.regnault@univ-reims.fr}

\date{July 2015}

\begin{abstract}
 We derive independence tests by means of dependence measures thresholding in a semiparametric context. 
 Precisely, estimates of 
 $\varphi$-mutual 
 informations, associated to $\varphi$-divergences between a joint distribution and the product distribution
 of its margins, are derived through  the dual representation of $\varphi$-divergences. The asymptotic properties of 
 the proposed estimates are established, including consistency, asymptotic distributions and large deviations principle. 
 The obtained tests of independence are 
  compared via their relative asymptotic Bahadur efficiency and numerical simulations. It follows that the proposed 
  semiparametric 
 Kullback-Leibler Mutual information test is the optimal one.  On the other hand, the proposed approach 
 provides a new method for estimating the Kullback-Leibler mutual information in a semiparametric setting,
 as well as a model selection procedure in  large class of  dependency models including semiparametric copulas.\\
 
\end{abstract}

\maketitle

\tableofcontents

\section{Introduction and notations}
Measuring the dependence between random variables has been a central aim of probability theory since its earliest developments. Classical examples 
of dependence measures are correlation measures of Pearson, Kendall or Spearman. While the first one 
focuses on linear relationship between real random variables, the two second ones measure the monotonic relationship between variables taking 
values in ordered sets. 
Pure-independence measures, between variables $X$ and $Y$ taking values in general measurable spaces 
$(\mathcal{X}, \mathcal{A}_{\mathcal{X}})$ 
and $(\mathcal{Y}, \mathcal{A}_{\mathcal{Y}})$, can be defined by considering any divergence between the joint distribution $\P$ of $(X,Y)$ and the 
product distribution
of its margins $\P^{\perp} :=\P_1 \otimes \P_2$, where $\P_1$ and $\P_2$ are, respectively, the marginal distributions of $X$ and $Y$. 
The most outstanding and widely used example of such  dependence measures is the $\chi^2$-divergence between $\P$ and $\P^\perp$ defined by
\begin{equation}\label{MIITEqnChiSq}
 \chi^2(\P, \P^\perp)   :=  \frac{1}{2} \int_{\mathcal{X}\times \mathcal{Y}} \left(\frac{\der \P}{\der \P^{\perp}}(x,y)-1 \right)^2 \der \P^\perp(x,y), 
\end{equation}
where $\frac{ \der\P}{\der \P^\perp}$ 
denotes the density of $\P$ with respect to  (w.r.t.) $\P^\perp$. Note that, if  $\P$ is a discrete distribution, i.e., if its support $\Om := \text{supp}(\P)$
 is discrete (finite or countably infinite) set,
 then the above divergence writes   
 $$\chi^2(\P, \P^\perp)   =   \frac{1}{2}\sum_{(x,y) \in \mathcal{X}\times \mathcal{Y}} \frac{\left(p_{x,y} - p_x p_y \right)^2}{p_x p_y}, $$
where $\P:=(p_{x,y})_{(x,y)}$, $\P^\perp =(p_x p_y)_{(x,y)}$, with $p_x :=\sum_{y} p_{x,y}$ and $p_y :=\sum_{x} p_{x,y}$.
Another classical example, associated to the Kullback-Leibler (KL) 
divergence between $\P$ and $\P^\perp$, is the well-known mutual information (MI)   defined by (see e.g. \cite{Cover_Thomas2006_Wiley})
\begin{equation} \label{MIITEqnMIClassical}
 I_{KL}(\P)  :=  \K(\P, \P^\perp):=\int_{\mathcal{X}\times \mathcal{Y}} \frac{\der \P}{\der \P^{\perp}}(x,y)\, \log \frac{\der \P}{\der \P^\perp}(x,y)\,  
 \der \P^\perp(x,y),
\end{equation}
which, in the case of discrete distributions, can be written under the form
$$I_{KL}(\P)=\sum_{(x,y)\in\mathcal{X}\times\mathcal{Y}}   p_{x,y}  \log \frac{p_{x,y}}{p_xp_y}.$$
We will call the above classical measures of dependence (\ref{MIITEqnChiSq}) and (\ref{MIITEqnMIClassical}), 
respectively,  $\chi^2$-mutual information ($\chi^2$-MI) 
and KL-mutual information (KL-MI).
When dealing with i.i.d. observations $(X_1,Y_1),\ldots, (X_n,Y_n)$,   of two random variables $(X,Y)$, we may 
test the null hypothesis, that the variables $X$ and $Y$ are independent, 
by means of estimating such dependence measure and deciding to reject the null 
hypothesis of independence if the estimate is sufficiently far from zero;
the classical $\chi^2$-independence test is such a procedure : the corresponding test statistic (in the discrete-distribution case) is
\begin{equation} \label{MIITEqnEstiChiEmp}
  2n\,  \chi^2\left( \widehat{\P}, \widehat{\P}^\perp \right) = n \sum_{(x,y)\in\Om} 
  \frac{\left(\widehat{p}_{x,y}- \widehat{p}_x \widehat{p}_y \right)^2}{ \widehat{p}_x 
  \widehat{p}_y},
\end{equation}
where $\widehat{\P} := \left(\widehat{p}_{x,y}\right)_{(x,y)}$ and $ \widehat{\P}^\perp :=\left(\widehat{p}_{x}\,\widehat{p}_{y}\right)_{(x,y)}$
are, respectively, the empirical versions of  $\P =(p_{x,y})_{(x,y)}$ and $\P^\perp = (p_xp_y)_{(x,y)}$.
Likewise, to test the independence, we can consider  as dependence measure the KL-MI and use the test statistic 
\begin{equation} \label{MIITEqnMIClassical_emp}
  2n \, I_{KL}(\widehat{\P}) = 
  2n\, \sum_{(x,y)\in\mathcal{X}\times\mathcal{Y}}   \widehat{p}_{x,y}  \log \frac{\widehat{p}_{x,y}}{\widehat{p}_x\widehat{p}_y}.
\end{equation}
The dependence measure can also be any other $\varphi$-divergence between $\P$ and $\P^\perp$. 
The tests based on such dependence measures, including the $\chi^2$-MI and KL-MI ones, have been extensively studied in the case 
of finite-discrete distributions; see e.g. \cite{Pardo2006_CH} Chapter 8, and the references therein.  
When dealing with continuous distributions  (or continuous random variables), 
obviously, the above direct plug-in estimates (\ref{MIITEqnEstiChiEmp}) and (\ref{MIITEqnMIClassical_emp}), of the dependence 
measures (\ref{MIITEqnChiSq}) and (\ref{MIITEqnMIClassical}), are not well defined.
Moreover, for countably-infinite discrete distributions, although the above estimates (\ref{MIITEqnEstiChiEmp}) and (\ref{MIITEqnMIClassical_emp}) 
remain well defined, their limiting distributions are not accessible.  Therefore, in the case of non finite-discrete distributions, particularly, for the widely used KL-MI,
other kind of estimates 
have been proposed and studied  in the literature; see e.g. \cite{Moon_Rajagopalan_Lall1995_PhysRevE} for a kernel density estimate, 
\cite{Kraskov_Stogbauer_Grassberger2004_PhysRevE} 
for a $k$-nearest-neighbor estimate extending those of Shannon entropy in one dimension 
based on $m$-spacing; see e.g. \cite{Tsybakov_Meulen1996_SJS}, \cite{Dudewicz_vanderMeulen1981_JASA} and
\cite{Beirlant_Dudewicz_Gyorfi_vanderMeulen1997_IJMSS} among others. 
\cite{Hulle_Marc2005_NC} derive an estimate using Edgeworth approximation of Shannon entropy. \cite{Darbellay_Vajda1999_IEEETIT},
\cite{Wang_Kulkarni_Verdu2005_IEEETIT} and   
\cite{Cellucci_Albano_Rapp2005_PhysRevE} propose estimates based on adaptative partitioning of $\mathcal{X}\times\mathcal{Y}$. See  also
\cite{Khan2007_PhysRevE} for an overview and numerical comparisons of these estimates. 
Based on the Kullback-Leibler importance estimation procedure,
 see \cite{Sugiyama_et_al2008_AISM},  \cite{Suzuki_Sugiyama_Sese_Kanamori2008} obtain
 an estimate of KL-MI called maximum likelihood mutual information, see also  
\cite{Sugiyama_Suzuki_Kanamori2012_CUP} Chapter 11.
Unfortunately, their (asymptotic) distributions remain 
inaccessible. Hence, testing independence from these estimates requires Monte-Carlo or Bootstrap approximations of the related $p$-values.
On the other hand, the above nonparametric estimates suffer from loss of efficiency, due to smoothing or partitioning, 
and suffer also from the difficulty of conveniently choosing the classes, the number of classes
or the  smoothing parameters (the bandwidths and the kernels).
The present paper introduces  new efficient semiparametric estimates of $\varphi$-mutual information ($\varphi$-MI), i.e., dependence measures 
associated to $\varphi$-divergence functionals, including the well known KL-MI and $\chi^2$-MI. 
These estimates are  obtained by making use of a dual representation of $\varphi$-MI, presented in 
Section~\ref{SecDuality}, without using any smoothing nor partitioning.  The obtained estimates are  
defined in the same way for both finite-discrete or non-discrete distributions, 
and coincide with the direct plug-in ones in the case of finite-discrete distributions. 
Their asymptotic properties are presented in Section~\ref{SecAsympProp}. Particularly, the consistency is stated for a large variety of 
semiparametric  models for $\der \P/\der \P^\perp$;
the asymptotic distribution is obtained for the KL-MI estimate in a 
special setting. 
The present approach  leads to new independence tests, 
whose Bahadur efficiency are compared in Section~\ref{SecBahadurEff} ; the most efficient test is shown to 
be the one based on the proposed estimate of  the particular KL-MI criterion.
It can be used also in order to build a large variety of dependence models, through for instance a cross validation-type 
model selection procedure based on the proposed estimate of $\varphi$-MI measure of dependence;
see Section \ref{Sec model select}.
 The powers of $\varphi$-MI based tests are compared numerically to classical noncorrelation tests in Section~\ref{SecSimus}.
 The results in the present paper have the advantage (unlike the classical noncorrelation tests) to remain valid in the case
of multisample problem (estimating $\varphi$-mutual informations of a multidimensional random variable
 as well as testing simultaneous independence of its components), but for simplicity, the results will be presented 
only for the two-sample case. The same results hold for the multisample problem. 
All proofs are postponed to the Appendix.

\section{$\varphi$-mutual informations, Dual representations and Estimation strategy} \label{SecDuality}
Given an i.i.d. sample, $(X_1,Y_1),\ldots, (X_n,Y_n)$, of a random vector $(X,Y)$ taking values in a measurable space $\left(\mathcal{X}\times
\mathcal{Y}, \mathcal{A}_{\mathcal{X}}\otimes
\mathcal{A}_{\mathcal{Y}} \right)$, we aim at testing the null hypothesis $\mathcal{H}_0$ of  independence of the margins $X$ and $Y$; formally
\begin{equation}\label{le probleme de test}
  \mathcal{H}_0 \,:\, X \textrm{ and } Y \textrm{ are independent}, \,\, \text{ against }\,\,  \mathcal{H}_1 \,:\,  X \textrm{ and } Y \textrm{ are dependent.}
\end{equation}
We derive such tests by estimating and thresholding $\varphi$-mutual informations between $X$ and $Y$ in a semiparametric context. Sections
\ref{SecDefMI}, \ref{SecModel} and \ref{SecEstStrat} to follow, respectively, define $\varphi$-mutual informations, present the semiparametric model 
under study, and introduce estimates of $\varphi$-MI used as test statistics for the test problem~(\ref{le probleme de test}).
Section \ref{Sec model select} defines a cross-validation procedure  for model selection
among $L$ candidate models for the ratio $\der \P/\der \P^\perp$, using the proposed estimate of $\varphi$-MI. 

\subsection{Introducing $\varphi$-mutual informations}\label{SecDefMI}
Denote by $\cM_1(\Om)$ the set of all probability distributions on the product measurable space 
$\left(\mathcal{X}\times\mathcal{Y}, \mathcal{A}_{\mathcal{X}}\otimes \mathcal{A}_{\mathcal{Y}} \right)$.
Let $\varphi : \R \to [0,+\infty]$ be some nonnegative closed proper convex function such that its domain 
$\text{dom}_\varphi:=\left\{x\in\R; \varphi(x)<\infty\right\}=:(a_\varphi,b_\varphi)$ is an interval, with endpoints $a_\varphi<1<b_\varphi$, and
$\varphi(1)=0$. The interval $(a_\varphi,b_\varphi)$ may be bounded or unbounded, open or not. 
The $\varphi$-divergence between any probability distributions $Q, P \in \cM_1(\Om)$, if $Q$ is absolutely continuous with respect to (a.c.w.r.t.) $P$, 
is defined by 
$$ D_\varphi(Q,P) := \int_{\Om}\varphi\left(\frac{d Q}{d P}(x,y)\right)\, d P(x,y).$$
If $Q$ is not a.c.w.r.t. $P$, we set $D_\varphi(Q,P)=+\infty$. 
Note that $D_\varphi(Q,P) \geq 0$, for any $Q$ and $P$. Moreover, if $\varphi$ is strictly convex on some neighborhood of $1$, we have the 
fundamental property
$$D_\varphi(Q,P) \geq 0, \text{ with equality if and only if }  Q=P.$$
In the following, we assume that the function $\varphi$ is  strictly convex and two times continuously differentiable on the interior of its domain 
$(a_\varphi, b_\varphi)$.  We have then $\varphi'(1)=0$, and without loss of generality, we can assume that $\varphi''(1)=1$.
The well-known Kullback-Leibler divergence $\K(\cdot,\cdot)$ is 
obtained for $\varphi(x) = \varphi_1(x):= x\log x -x+1$, the ``modified'' Kullback-Leibler
divergence $\K_m(\cdot,\cdot)$ is obtained for $\varphi(x)=\varphi_0(x):=-\log x +x -1$. The $\chi^2$ and modified-$\chi^2$ divergences, denoted 
$\chi^2(\cdot,\cdot)$
and $\chi^2_m(\cdot,\cdot)$, are associated, respectively, to the convex functions $\varphi(x)=\varphi_2(x):=(x-1)^2/2$ and 
$\varphi(x) = \varphi_{-1}(x) := (x-1)^2/(2x)$. The so-called  
Hellinger distance $H(\cdot,\cdot)$ is obtained for $\varphi(x)= \varphi_{1/2}(x) :=2 (\sqrt{x}-1)^2$; see Table~\ref{table convex conjugates}.
All  these divergences are members of the so-called ``power-divergences''
$D_{\varphi_\gamma}(\cdot,\cdot)$ associated to the convex functions  $\varphi_\gamma(\cdot)$ defined by
\begin{equation} \label{MITEqnPowerFunction}
 \varphi_\gamma(\cdot) : x\in\mathbb{R}_{+}^{\ast}\mapsto\varphi_{\gamma}(x):=\frac{x^{\gamma}-\gamma
x+\gamma-1}{\gamma(\gamma-1)} 
\end{equation}
if $\gamma\in\mathbb{R}\setminus\left\{  0,1\right\}  $, $\varphi
_{0}(x):=-\log x+x-1$ and $\varphi_{1}(x):=x\log x-x+1$.
The standard divergences $\K(\cdot,\cdot)$,  $\K_m(\cdot,\cdot)$, $\chi^2(\cdot,\cdot)$,
 $\chi^2_m(\cdot,\cdot)$ and $H(\cdot,\cdot)$ are then associated, respectively, to the real convex functions
 $\varphi_1(\cdot)$, $\varphi_0(\cdot)$, $\varphi_2(\cdot)$, $\varphi_{-1}(\cdot)$ and $\varphi_{1/2}(\cdot)$.
Note that the divergences are generally not symmetric; particularly, we have for any $Q, P \in \cM_1(\Om)$, 
 $\K_m(Q,P) = \K(P,Q)$ and $\chi^2_m(Q,P) = \chi^2(P,Q)$.
For more details and proofs, we can refer to \cite{Liese_Vajda1987_BSB} and \cite{Broniatowski_Keziou2006_STUDIA}.
For any probability distribution 
$P \in \cM_1(\Om)$, let $P^\perp$ denotes the product distribution $P^\perp:=P_1\otimes P_2$ of the margins 
$P_1$ and $P_2$ of $P$. 
The $\varphi$-mutual information of $P$, associated to the divergence $D_\varphi(\cdot,\cdot)$, is defined as
$$ I_{\varphi}(P) := D_{\varphi} (P, P^{\perp}).$$
For any random vector $(X,Y)$  defined on a probability space $(\Omega,\mathcal{A},\mathbf{P})$ and taking its values in $\left(\mathcal{X}\times
\mathcal{Y}, \mathcal{A}_{\mathcal{X}}\otimes
\mathcal{A}_{\mathcal{Y}} \right)$, with joint distribution $\P\in \cM_1\left(\mathcal{X}\times\mathcal{Y}\right)$, 
the $\varphi$-mutual information ($\varphi$-MI) of $(X,Y)$ is defined to be
\begin{equation}\label{phi-MI de (X,Y)}
I_\varphi(X,Y) := I_\varphi(\P) = D_{\varphi} (\P, \P^{\perp}) = \int_{\Om}\varphi\left(\frac{\der \P}{\der \P^\perp}(x,y)\right)\, \der \P^\perp(x,y).
\end{equation}
Since $D_\varphi(\P,\P^\perp)\geq 0$, with equality if and only if $\P=\P^\perp$, i.e., if and only if $X$ and $Y$ are independent, $\varphi$-MI measures 
then the dependence between the random variables $X$ and $Y$. In contrast to the correlation coefficients of Pearson, Kendall or Spearman, the 
$\varphi$-MI does not focus on the linear or monotonic relationship between random variables; it constitutes a proper dependency measure. 
Note that $I_{\varphi_1}$ and $I_{\varphi_2}$, with $\varphi_1$ and  $\varphi_2$ given in Table~\ref{table convex conjugates}, are, respectively, the  KL-MI   
and $\chi^2$-MI, given by (\ref{MIITEqnMIClassical}) and~(\ref{MIITEqnChiSq}). 
Thus, the test problem (\ref{le probleme de test}) is equivalent, in the context of $I_\varphi$ criteria, to testing
\begin{equation*}
   I_\varphi(\P) =0 \quad \text{ against }\quad    I_\varphi(\P) > 0.
\end{equation*}
Hence, we can use as test statistic an estimate of $I_{\varphi}(\P)$, and reject the null hypothesis $\mathcal{H}_0$ when the estimate takes large 
values. 
A natural attempt to estimate the $\varphi$-MI of $(X,Y)$ consists in considering the plug-in estimate of $I_\varphi(\P)$ obtained by replacing 
$\P(\cdot)$ by its empirical counterpart 
\begin{equation} \label{MITEqnProbEmp}
 \widehat{\P}(\cdot) = \frac{1}{n} \sum_{i=1}^n \delta_{(X_i,Y_i)}(\cdot),
\end{equation}
associated to the i.i.d. sample $(X_1,Y_1), \ldots, (X_n,Y_n)$ of $(X,Y)$. Here, $ \delta_{(x,y)}(\cdot)$ denotes the Dirac measure
at $(x,y)$ for all $(x,y) \in \Om$.
Unfortunately, by doing so, we only measure dependence of the contingency table associated to the sample. When dealing with variables $X$ and $Y$ 
absolutely continuous with respect to Lebesgue measure, the contingency table is almost surely an $n\times n$ table with all coefficients 
except diagonal 
ones equal to zero ; particularly, variables $X$ and $Y$ appear (misleadingly) purely dependent, yielding to reject systematically the null hypothesis.
A second, less crude, approach consists in gathering the values $X_i$ and $Y_i$ into classes
and testing independence between the induced finite-discrete variables $\tilde{X}$ and $\tilde{Y}$, by empirically estimating the $\varphi$-MI of 
$(\tilde{X}, \tilde{Y})$. This widespread approach  suffers from the difficulty of conveniently choosing the classes. Moreover, an important amount of 
information carried by the sample is lost during this process, yielding to poor efficiency -- or power -- of these tests. An other approach, is to use 
 kernel nonparametric estimates of the joint density and the marginal ones, but as it is well known this provides 
less efficient estimates and leads to the difficulty of choosing the optimal smoothing parameters. 
As an alternative, we propose in the present paper semiparametric modeling of the ratio $\der\P/\der \P^\perp$, and the use of duality
to obtain well-defined estimates of $\varphi$-MI without smoothing nor partitioning. The present approach applies 
for both continuous or discrete distributions, or mixtures of continuous and discrete distirubtions. 

\subsection{Semiparametric modeling of the ratio $\der\P/\der\P^\perp$}\label{SecModel}
Assume that the joint distribution $\P$ of the random vector $(X,Y)$ belongs to the semiparametric model 
\begin{equation}\label{le modele}
 \cM_{\Theta} := \left\{P \in \cM_1(\Om) \text{ such that } 
 \frac{d P}{d P^{\perp}}(\cdot,\cdot) = : h_\theta(\cdot,\cdot);\,  \theta \in \Theta \right\},
\end{equation}
where $\Theta \subset \R^{1+d}$ is the  parameter space, and $h_{\theta}(\cdot,\cdot) : (x,y)\in\Om \mapsto h_\theta(x,y)\in\R$ is some specified 
real-valued function, indexed by the parameter $\theta$. 
In the sequel, we will consider the following assumptions on the model $\cM_{\Theta}$.
\begin{itemize}
 \item[(A.1)] $\left( h_{\theta}(x,y) = h_{\theta'}(x,y), \forall (x,y) \in \Om \right) \, \, \Rightarrow \, \,  \left( \theta=\theta' \right)$ (identifiability); 
 \item[(A.2)] there exists (a unique) $\theta_0 \in \text{int}(\Theta)$ satisfying $h_{\theta_0}(x,y) = 1$, $\forall (x,y) \in \mathcal{X} \times \mathcal{Y}$.
\end{itemize}

Assumption (A.1) is a natural identifiability condition for $d P/{d P}^\perp$. 
Assumption (A.2) ensures independence is covered by  
the model $\cM_{\Theta}$. The uniqueness of $\theta_0$ follows from Assumption (A.1). 
Denote by $\theta_T$ the ``true'' unknown value of the parameter, namely, the unique value satisfying 
$$ \frac{\der\P}{\der\P^\perp} (x,y) = h_{\theta_T}(x,y),\quad \forall (x,y)\in\mathcal{X}\times\mathcal{Y},$$
which is assumed to be an interior point of $\Theta$.  Then, we have $\theta_T=\theta_0$ if and only if   $X$
and $Y$ are independent. 
Below are listed some relevant examples of the model $(\ref{le modele})$.
\begin{example}  \label{modele gaussien}
Let $(X,Y)\in\mathbb{R}^2$ be a centered Gaussian random vector with correlation coefficient $\rho \in]-1,1[$
and centered normal margins with the same variance $\sigma^2>0$.  
A straightforward computation shows that the ratio $\mathrm{d} \P / \mathrm{d} \P^\perp$ can be written under the form of the model (\ref{le modele})
where
\begin{equation}
h_\theta(x,y) = \exp\left\{\alpha+\beta_1(x^2+y^2)+\beta_2xy\right\},
\end{equation}
$\theta:=(\alpha,\beta_1,\beta_2)^\top \in \R^3$, with $\alpha=-\log({1-\rho^2})/2,$  $\beta_1 = -\rho^2/(2\sigma^2(1-\rho^2))$  and $\beta_2 = \rho/
(\sigma^2(1-\rho^2))$. Note that the parameter value,  corresponding to the independence hypothesis, is  $\theta_0=(0,0,0)^\top.$
Moreover, if the distribution of $(X,Y)$ is Gaussian with unknown mean $\mu:=(\mu_1,\mu_2)\top$ and unknown variance 
matrix $\Gamma$, then we can 
show that the ratio $\mathrm{d} \P/ \mathrm{d} \P^\perp$ can be written under the form of the model (\ref{le modele}) with
\begin{equation}
h_\theta(x,y) = \exp\left\{\alpha+\beta_1x+\beta_2y+\beta_3x^2+\beta_4y^2+\beta_5xy\right\},
\end{equation}
and $\theta := (\alpha,\beta_1,\beta_2,\beta_3,\beta_4,\beta_5)^\top$. 
Note that the number of free 
parameters  in $\theta_T$ is $d=5$, and that $\alpha_T$ is considered as a normalizing parameter due to the constraint 
$\int_{\Om} h_{\theta_T}(x,y)\, \mathrm{d} \P^\perp(x,y)=\int_{\Om} \mathrm{d} \P(x,y)=1$ since $\P$ is a probability distribution.  
Moreover, we have $\theta_0 =(0,\ldots,0)^\top\in\R^6.$
 \end{example} 

\begin{example} \label{Le modele general}Let $\psi_0(\cdot,\cdot):=\mathds{1}_{\Om}(\cdot,\cdot), \psi_1(\cdot,\cdot),  \psi_2(\cdot,\cdot), \ldots,$ be 
some basis functions of the space $L^2(\mathcal{X}\times\mathcal{Y},\mathbb{P}^\perp)$, 
and assume that $\log({\mathrm{d}\P}/{\mathrm{d}\P^\perp}(\cdot,\cdot))\in 
L^2\left(\Om,\P^\perp\right)$. 
We can then build increasing models of the form (\ref{le modele}) developing the function 
$$(x,y)\in\mathcal{X}\times\mathcal{Y}\mapsto \log \frac{\mathrm{d}\mathbb{P}}{\mathrm{d}\mathbb{P}^\perp}(x,y)$$
according to the above basis functions. Using for instance the first $(1+d)$-basis functions, we obtain the following model for  
${\mathrm{d} \mathbb{P}} / {\mathrm{d} \mathbb{P}^\perp}(\cdot,\cdot)$
$$h_\theta : (x,y)\in\Om\mapsto h_\theta(x,y)=\exp\left(\alpha+\beta_1\psi_1(x,y)+\cdots +\beta_d\psi_d(x,y)\right),$$
\end{example} 
where $\theta = 
(\alpha,\beta_1,\ldots,\beta_d)^\top\in\Theta\subset{\R}^{1+d}$. Then, the independence parameter value is $\theta_0=(0,\ldots,0)^\top\in\R^{1+d}.$

\begin{example} \label{modele fini}
 Assume that the support of $\P$, $\text{supp}(\P) =: \mathcal{X}\times\mathcal{Y}$, is a known finite-discrete 
 set of size $K_1K_2$; denote by $\left(\P(x,y)\right)_{(x,y)\in\mathcal{X}\times
 \mathcal{Y}} :=
 \left(p_{x,y}\right)_{(x,y)\in\mathcal{X}\times\mathcal{Y}}$ the density of $\P$ with respect to the counting measure on $\Om$. Then we have 
 \begin{equation}\label{modele exponentiel fini}
\frac{\mathrm{d} \P}{\mathrm{d} \P^{\perp}}(x,y) = \exp\left( \sum_{(a,b) \in \Om} \theta_{a,b}\,\,  \mathds{1}_{\{a\}}(x) \, \mathds{1}_{\{b\}}(y) \right),
\end{equation}
 where 
 $$\theta_{a,b} = \log \frac{p_{a,b}}{p_a p_b}, \quad (a,b) \in \Om.$$
 If we denote for instance the elements of $\mathcal{X}$ and $\mathcal{Y}$ as follows
 $$\mathcal{X}:=\left\{a_1,\ldots,a_{K_1}\right\} \quad\text{and}\quad \mathcal{Y}:=\left\{b_1,\ldots,b_{K_2}\right\},$$ 
then we can see that $\P$ belongs to the model (\ref{le modele})
taking 
\begin{equation}\label{modele exponentiel fini bis}
h_{\theta}(x,y)=\exp\left(\alpha + \sum_{(i,j)\neq (1,1)} \beta_{i,j} \, \mathds{1}_{\{a_i\}}(x)\, \mathds{1}_{\{b_j\}}(y)\right),
\end{equation}
with the parametrization  $\theta = (\alpha,\beta^\top)^\top \in \R^{K_1K_2}$, where $\alpha$ is a scalar and $\beta = (\beta_{i,j})_{(i,j)\neq (1,1)}$ is the 
$(K_1K_2-1)$-dimensional vector obtained from the $K_1\times K_2$-matrix of real entries $(\beta_{i,j})$ removing the first entry 
$\beta_{1,1}$. Moreover, we have for the true value $\theta_T$ 
$$\alpha_T = \log  \frac{p_{a_1,b_1}}{p_{a_1} p_{b_1}}, \quad \text{and}\quad {\beta_{i,j}}_T = \log \frac{p_{a_i,b_j}}{p_{a_i} p_{b_j}} - 
\log\frac{p_{a_1,b_1}}{p_{a_1} p_{b_1}},$$
for all $(i,j)\in\{1,\ldots,K_1\}\times\{1,\ldots,K_2\}\setminus \{(1,1)\}$, and that
 the number of free-parameters in $\theta_T$ is 
 equal to $(K_1-1)(K_2-1)$. Moreover, we have $\theta_0=(0,\ldots,0)^\top\in \R^{K_1K_2}$.
\end{example}

\begin{example} Assume that the distribution $\P$ of the random vector $(X,Y)\in\R^2$ is 
of continuous margins.
The copula $C(\cdot,\cdot)$ of the vector $(X,Y)$, see e.g. 
\cite{Nelsen2006_Springer},
is defined,  $\forall (u,v)\in ]0,1[^2$, by 
$$C(u,v) := F(F_1^{-1}(u), F_2^{-1}(v)),$$
where $F(\cdot,\cdot)$ is the cumulative distribution function of the vector $(X,Y)$, and $F_1$ and $F_2$ are the (marginal)
cumulative distribution functions of $X$ and $Y$, respectively.  
The copula $C(\cdot,\cdot)$ is in itself a distribution function on  $ ]0,1[^2$. If 
$F(\cdot,\cdot)$ is absolutely continuous with respect to the Lebesgue measure on $\R^2$, 
then we have the relation
$$ \frac{\mathrm{d} \P}{\mathrm{d} \P^\perp}(x,y)=\frac{f(x,y)}{f_1(x)f_2(y)} = {c\left(F_1(x),F_2(y)\right)},$$
where $f(\cdot,\cdot)$ is the joint density of $(X,Y)$,  $f_1$ and  $f_2$ are  the marginal densities of $X$ and $Y$, 
and $c(\cdot,\cdot)$ the copula density. 
Numerous parametric examples of the model (\ref{le modele}) can then be obtained taking the function
\begin{equation}\label{examples copules}
h_\theta(x,y) = {c_\beta(F_{1,\gamma_1}(x),F_{2,\gamma_2}(y))}
\end{equation}
where $\left\{c_\beta(\cdot,\cdot);\, \beta\in D\subset\R^m\right\}$ is some parametric copula density model, see e.g. 
\cite{Nelsen2006_Springer} or 
\cite{Joe1997_Chapman}
for examples of such models, and $\left\{F_{1,\gamma_1}; \gamma_1\in\Gamma_1\right\}$ and $\left\{F_{2,\gamma_2}; \gamma_2\in
\Gamma_2\right\}$
are some parametric models for the marginal distribution functions. Here, the parameter  of interest is $\theta := 
(\gamma_1,\gamma_2,\beta)\in
\Theta :=\Gamma_1\times\Gamma_2\times D$. Note that the assumption (A.2) is generally not satisfied for this particular model. 
In fact, if we denote $\beta_0$ 
the particular value corresponding to the copula of independence, then we have $h_{(\gamma_1,\gamma_2,\beta_0)}(\cdot,\cdot) =1$
for any $(\gamma_1,\gamma_2)\in\Gamma_1\times\Gamma_2$. Although assumption (A.2) is generally not satisfied, models 
(\ref{examples copules}) can be used in estimating $\varphi$-MI under the assumption that the margins are dependent.  
\end{example} 

\begin{example} \label{modele copule nonparam} We can also deal with  semiparametric models induced by semiparametric models of 
copula densities,
with nonparametric unknown continuous  marginal distribution functions $F_1(\cdot)$ and $F_2(\cdot)$, taking
$$h_{\theta}(x,y) = {c_\theta(F_1(x),F_2(y))};\, \theta\in\Theta\subset\R^d.$$
\end{example} 

\subsection{Dual representation and dual estimation of $\varphi$-MI} \label{SecEstStrat}
We define estimates of $\varphi$-MI by taking advantage of the modeling~(\ref{le modele}) and the dual representation  of $\varphi$-divergences 
obtained in 
 \cite{Keziou2003_CRAS} and \cite{Broniatowski_Keziou2006_STUDIA}.
Denote $\varphi^*(\cdot)$ the convex conjugate of the convex function $\varphi(\cdot)$, namely, the function defined by 
$$\varphi^* : t\in\R\mapsto \varphi^*(t):=\sup_{x\in\R}\left\{tx-\varphi(x)\right\}\in\R\cup\{+\infty\}.$$
Note that $\varphi^*(\cdot)$ is, in turn, a proper closed convex function, in particular, $\varphi^*(0)=0.$ 
Assume that $\varphi(\cdot)$ is essentially smooth, i.e., differentiable on $]a_\varphi,b_\varphi[$ with 
$\lim_{x\downarrow a_\varphi}\varphi'(x) = -\infty$ if $a_\varphi>-\infty$ and 
$\lim_{x\uparrow b_\varphi}\varphi'(x) = +\infty$ if $b_\varphi<+\infty$. This is equivalent to the condition that 
$\varphi^*(\cdot)$ is strictly convex on its domain. Provided that 
\begin{itemize}
\item [(A.3)] the $\varphi$-mutual information  $I_\varphi(\P)<\infty$, 
\end{itemize}
see its definition (\ref{phi-MI de (X,Y)}), it can be rewritten under the form
\begin{equation}\label{eqn 1}
I_\varphi(\P) = \sup_{f \in \mathcal{F}}\left\{\int_{\Om} f(x,y) \, \der\P(x,y) - \int_{\Om} \varphi^*\left(f(x,y)\right)\, \der\P^\perp(x,y)\right\},
\end{equation}
where $\mathcal{F}$ is any class, of measurable real-valued functions  $f : \Om \to \overline{\R}$, 
that contains the particular function $\varphi'(\der\P/\der\P^\perp)$ and 
satisfies  the condition $\int_{\Om} |f|\, \der\P <\infty$,
for all $f\in\mathcal{F}$.   
Note that,  for all $x \in (a_{\varphi}, b_{\varphi})$, we have
$$ \varphi^*(\varphi'(x)) = x \varphi'(x)-\varphi(x).$$
In Table \ref{table convex conjugates} are given explicit formulas of convex conjugates of some standard divergences.
\begin{table}[h] 
 \centering
 \begin{tabular}{|l||l|l||l|l|} 
   \hline
  $D_\varphi(\cdot,\cdot)$ & $\varphi(\cdot)$ & $\text{dom}\varphi $ & $\text{dom}\varphi^* $ & $\varphi^*(\cdot)$ \\
  \hline
  \hline
  $\K_m(\cdot,\cdot)$ &  $\varphi_0(x):=-\log x +x -1$ & $]0,+\infty[$  & $]-\infty,1[$ & $ - \log(1-t)$ \\
  \hline
  $\K(\cdot,\cdot)$ &  $\varphi_1(x):=x\log x -x +1$ & $[0,+\infty[$ & $\mathbb{R}$ & $e^t-1$  \\
  \hline
  $\chi^2_m(\cdot,\cdot)$ &  $\varphi_{-1}(x):=\frac{1}{2}\frac{\left(x-1\right)^2}{x}$
  & $]0,+\infty[$ &
  $\left]-\infty,\frac{1}{2}\right]$ & $1-\sqrt{1-2t}$ \\
  \hline
  $\chi^2(\cdot,\cdot)$ &  $\varphi_2(x):=\frac{1}{2}\left(x-1\right)^2$ &  $\mathbb{R}$ &  $\mathbb{R}$ & $\frac{1}{2}t^2+t$ \\
  \hline
  $H(\cdot,\cdot)$ &  $\varphi_{1/2}(x):=2(\sqrt{x}-1)^2$ & $[0,+\infty[$ &  $]-\infty, 2[$ & $\frac{2t}{2-t}$ \\
  \hline
 \end{tabular}
 \vskip 0.2cm
 \caption{Convex conjugates for some standard divergences.}
   \label{table convex conjugates} 
\end{table} 
From (\ref{eqn 1}), taking into account the model (\ref{le modele}) by specifying 
$$\mathcal{F}=\{\varphi'(h_{\theta});\,  \theta \in \Theta\},$$
and assuming in addition that
\begin{itemize}
 \item[(A.4)] for all $\theta \in \Theta$, we have $\int_{\Om} \left|\varphi'(h_{\theta}(x,y))\right| \, \der\P(x,y) < \infty$,
\end{itemize}
we obtain
\begin{equation}\label{formule 1}
 I_\varphi(\P) = \sup_{\theta\in\Theta}\left\{\int_{\Om}\varphi'(h_\theta(x,y))\, \der\P(x,y) - \int_{\Om}\varphi^*\left(\varphi'(h_\theta(x,y))\right)\, \der\P^\perp(x,y)
 \right\}.
\end{equation}
Moreover, the supremum is unique and achieved in $\theta=\theta_T$. The uniqueness of the supremum $\theta_T$ follows from
the strict convexity of $\varphi^*(\cdot)$ and the identifiability assumption (A.1).
We propose then the following ``dual'' estimate of $I_\varphi(\P)$
\begin{eqnarray}\label{estim div}
\widehat{I}_\varphi & := & \sup_{\theta\in\Theta} \left\{\int_{\Om}\varphi'\left(h_\theta(x,y)\right)\, \der\widehat{\P}(x,y) - \int_{\Om}
\varphi^*\left(\varphi'(h_\theta(x,y))\right)\, \der\widehat{\P}_1\otimes\widehat{\P}_2(x,y)\right\}\nonumber\\
  & = & \sup_{\theta\in\Theta} \left\{\frac{1}{n}\sum_{i=1}^n \varphi'\left(h_\theta(X_i,Y_i)\right) - \frac{1}{n^2}\sum_{i=1}^n\sum_{j=1}^n
  \varphi^*\left(\varphi'\left(h_\theta(X_i,Y_j)\right)\right)\right\},
\end{eqnarray}
and the following ``dual'' estimate of the parameter $\theta_T$
\begin{eqnarray}\label{estim param}
\widehat{\theta}_\varphi & := & \arg\sup_{\theta\in\Theta} \left\{\int_{\Om}\varphi'\left(h_\theta(x,y)\right)\, \der\widehat{\P}(x,y) - \int_{\Om}
\varphi^*\left(\varphi'(h_\theta(x,y))\right)\, \der\widehat{\P}_1\otimes\widehat{\P}_2(x,y)\right\}\nonumber\\
 & = & \arg\sup_{\theta\in\Theta} \left\{\frac{1}{n}\sum_{i=1}^n \varphi'\left(h_\theta(X_i,Y_i)\right) - \frac{1}{n^2}\sum_{i=1}^n\sum_{j=1}^n
  \varphi^*\left(\varphi'\left(h_\theta(X_i,Y_j)\right)\right)\right\}, 
\end{eqnarray}
where $\widehat{\P}(\cdot)$ is the empirical distribution, associated to the sample, given by (\ref{MITEqnProbEmp}).
For ease of presentation, 
define, $\forall\theta\in\Theta$ and $\forall(x,y)\in\Om$, the functions 
\begin{equation}\label{la fonction f}
f_\theta(x,y) := \varphi'(h_\theta(x,y)),
\end{equation}
\begin{equation}\label{la fonction g}
g_\theta(x,y) := \varphi^*\left(\varphi'(h_\theta(x,y))\right) = h_\theta(x,y) \varphi'\left(h_\theta(x,y)
\right) - \varphi\left(h_\theta(x,y\right)),
\end{equation}
 which we assume to be continuous, in $\theta$, on the set $\Theta$,  
\begin{equation}\label{la fonction M theta}
 M : \theta\in\Theta\mapsto M(\theta) :=   \int_{\Om}f_\theta(x,y)\, \der\P(x,y) - \int_{\Om}
g_\theta(x,y)\, \der\P_1\otimes\P_2(x,y)
\end{equation}
and its empirical version 
\begin{equation}\label{la fonction Mn theta}
 M_n : \theta\in\Theta\mapsto M_n(\theta) :=   \int_{\Om}f_\theta(x,y)\, \der\widehat{\P}(x,y) - \int_{\Om}
g_\theta(x,y)\, \der\widehat{\P}_1\otimes\widehat{\P}_2(x,y).
\end{equation}
Therefore, the formula (\ref{formule 1}) becomes
 \begin{equation}\label{formule 1 bis}
 I_\varphi(\P) = \sup_{\theta\in\Theta} \, M(\theta) = M(\theta_T), \quad \text{and} \quad \theta_T = \arg\sup_{\theta\in\Theta}\, M(\theta).
\end{equation}
The estimates (\ref{estim div}) and (\ref{estim param}), in turn, can be written as
 \begin{equation}\label{estim div bis}
\widehat{I}_\varphi = \sup_{\theta\in\Theta} \, M_n(\theta) = M_n(\widehat{\theta}_\varphi) 
\end{equation}
and 
 \begin{equation}\label{estim param bis}
\widehat{\theta}_\varphi = \arg\sup_{\theta\in\Theta} \, M_n(\theta).
\end{equation}
Note that the functions $f_\theta(\cdot,\cdot)$, $g_\theta(\cdot,\cdot)$,  $M(\cdot)$ and $M_n(\cdot)$
all depend on $\varphi(\cdot)$, but the subscript $\varphi$ is omitted for simplicity.

\begin{example} \label{MITExpDualEstFinite}
 In the context of finite-discrete distributions, using the exponential model described in Example \ref{modele fini}, 
 we show that the proposed dual estimate (\ref{estim div}) of 
 $I_{\varphi}(\P)$, obtained by the above ``duality'' technique, equals the direct plug-in one 
 \begin{equation}\label{estimMIempDirect}
 \widehat{I}_{\varphi}^{\textrm{emp}} :=  I_\varphi(\widehat{\P}) = \sum_{(x,y) \in \Om} \varphi \left( 
 \frac{\widehat{p}_{x,y}}
 {\widehat{p}_x \widehat{p_y}} \right) 
 \widehat{p}_x \widehat{p}_y.
 \end{equation}
 Indeed, we have by its proper definition
 \begin{equation} \label{EqnScoreEmp}
 \widehat{I}_{\varphi} = \sup_{\theta \in \Theta} M_n(\theta), \text{ where } M_n(\theta) = \sum_{(x,y) \in \Om} \left[\varphi'(e^{\theta_{x,y}})
 \widehat{p}_{x,y} - e^{\theta_{x,y}}\varphi'(e^{\theta_{x,y}})\widehat{p}_x \widehat{p}_y 
  + \varphi(e^{\theta_{x,y}})\widehat{p}_x \widehat{p}_y \right].
 \end{equation}
 Differentiating~(\ref{EqnScoreEmp}) with respect to $\theta_{x,y}$ for $(x,y) \in \Om$ yields
 $$\frac{\partial}{\partial \theta_{x,y}} M_n(\theta) = \varphi''(e^{\theta_{x,y}}) \left( e^{\theta_{x,y}}\, \widehat{p}_{x,y}  -  
 e^{2\theta_{x,y}} \,\widehat{p}_x \widehat{p}_y \right).$$
 Canceling derivatives $\frac{\partial}{\partial \theta_{x,y}} M_n(\theta)$ yields
  $$\widehat{\theta}_{x,y}=\log \frac{\widehat{p}_{x,y}}{\widehat{p}_x \widehat{p}_y}, \quad (x,y) \in \Om,$$
  which is independent from the choice of $\varphi$ for this particular model.
 Finally, straightforward simplifications yield
  $$\widehat{I}_{\varphi}  =  M_n(\widehat{\theta}) 
    =  \sum_{(x,y) \in \Om} \varphi \left( \frac{\widehat{p}_{x,y}}{\widehat{p}_x \widehat{p}_y} \right) \widehat{p}_x \widehat{p}_y 
 =  \widehat{I}_{\varphi}^{\textrm{emp}}.$$
 Particularly, for $\varphi(x)=\varphi_2(x):=(x-1)^2/2$,  
 the estimate  $\widehat{I}_{\varphi_2}$  of the $\chi^2$-mutual information -- or $\chi^2$ measure of independence -- obtained by 
 the duality technique is shown to equal (up to the factor $2n$) the classical $\chi^2$ statistics. Hence, in the context of finite-discrete distributions,
 using the exponential model described  in Example \ref{modele fini}, we see that the proposed approach, via duality technique, 
 recovers the classical direct plug-in one, in particular, the well-known classical $\chi^2$-independence test.
\end{example}

\begin{remark}  For finite discrete distributions (with known support, of size say $K$, see Example \ref{modele fini}), 
as in plug-in estimation of Shannon entropy (see e.g. \cite{Chao_Shen2003_EES}), 
the direct plug-in estimates $\widehat{I}_\varphi^{emp}$ are valid with small bias if 
the sample size $n>>K$. If  the sample size $n$ is not sufficiently large compared to the space size $K$,  
models $h_\theta(\cdot)$ other than (\ref{modele exponentiel fini})  
should be used (through e.g. the model selection 
procedure described in Section \ref{Sec model select}), with small parameter 
dimension, 
and the corresponding
dual estimate $\widehat{I}_\varphi$, if the model $h_\theta(\cdot)$ is correctly specified,  
could be more promising than  the direct plug-in one 
$\widehat{I}_\varphi^{emp}$. 
\end{remark}

\begin{example}
Note that when dealing with semiparametric copula models 
$$h_{\theta}(x,y) = {c_{\theta} (F_1(x), F_2(y))},$$ 
with  unknown nonparametric cumulative distribution 
functions $F_1$ and $F_2$, it is necessary to estimate them, using for example their empirical counterparts.
Denote  by $\widehat{F}_1(\cdot)$ and  $\widehat{F}_2(\cdot)$ the empirical cumulative distribution functions associated, respectively, to   
the samples $X_1\ldots,X_n$ and $Y_1,\ldots,Y_n$, i.e.,
$$\widehat{F}_1(x) := \frac{1}{n} \sum_{i=1}^n \mathds{1}_{]-\infty,x]}(X_i) \quad \text{and} \quad 
\widehat{F}_2(y) := \frac{1}{n} \sum_{i=1}^n \mathds{1}_{]-\infty,y]}(Y_i).$$
So that $\widehat{I}_{\varphi}$ and $\widehat{\theta}_{\varphi}$ become
$$ \widehat{I}_{\varphi}  =  \sup_{\theta\in\Theta} 
 \left\{ \frac{1}{n} \sum_{i=1}^n \varphi'\left( {c_{\theta} \left( \widehat{F}_1(X_i),\widehat{F}_2(Y_i)\right)} \right) - \frac{1}
 {n^2} 
 \sum_{i=1}^n \sum_{j=1}^n  \varphi^* \left( \varphi' \left( {c_{\theta} \left( \widehat{F}_1(X_i),\widehat{F}_2(Y_j) \right)}\right) \right) \right\}$$
$$  \widehat{\theta}_{\varphi}  =  \arg\sup_{\theta\in\Theta} 
 \left\{ \frac{1}{n} \sum_{i=1}^n \varphi'\left( {c_{\theta} \left( \widehat{F}_1(X_i),\widehat{F}_2(Y_i)\right)} \right) - \frac{1}
 {n^2} 
 \sum_{i=1}^n \sum_{j=1}^n  \varphi^* \left( \varphi' \left( {c_{\theta} \left( \widehat{F}_1(X_i),\widehat{F}_2(Y_j) \right)}\right) \right) \right\}.$$
Note that $n\widehat{F}_1(X_i)$ is the rank of $X_i$ in the sample $X_1, \dots, X_n$ and $n\widehat{F}_2(X_j)$ 
is the rank of $Y_j$ in the sample 
$Y_1, \dots, Y_n$.  For some copula models, the copula density $c_\theta(u_1,u_2)$ may be unbounded 
when either $u_1$ or $u_2$ tends to 1; see e.g. 
\cite{Genest_Ghoudi_Rivest1995_Biometrika}. In this case, to avoid this difficulty, the ``rescaled'' empirical cumulative distribution functions
$$\widetilde{F}_1(\cdot) := \frac{n}{n+1} \widehat{F}_1(\cdot), \quad \widetilde{F}_2(\cdot) := \frac{n}{n+1} \widehat{F}_2(\cdot)$$
should be used instead of the standard ones $\widehat{F}_1(\cdot)$ and $\widehat{F}_2(\cdot)$.
\end{example}

\subsection{A model selection procedure for the ratio ${\der \P}/{\der \P^\perp}$  through $\varphi$-MI criterion}\label{Sec model select}
Let $\cM_{\Theta_1}:=\left\{ h_{\theta_1,1}(\cdot,\cdot); \, \theta_1\in\Theta_1\subset\R^{d_1} \right\}, \ldots, \cM_{\Theta_L}:=\left\{ h_{\theta_L,L}(\cdot,\cdot); 
\, \theta_L\in \Theta_L\subset\R^{d_L} \right\}$ be $L$ candidate models for the ratio ${\der \P}/{\der\P^\perp}$. 
For any model $\mathcal{M}_{\Theta_\ell}$, denote by $\widehat{\theta}_\ell$ the estimate of $\theta_T$ given by 
$$\widehat{\theta}_\ell :=\arg\sup_{\theta_\ell\in\Theta_\ell} M_n(\theta_\ell).$$
The corresponding ``expected'' criterion is 
$$M(\widehat{\theta}_\ell) = \int_{\Om} f_{\widehat{\theta}_\ell}(x,y)\, \der\P(x,y) - \int_{\Om} g_{\widehat{\theta}_\ell}(x,y)\, \der\P(x,y)^\perp.$$
 From the representation (\ref{formule 1 bis}), we can see that the larger the expected criterion $M(\widehat{\theta}_\ell)$
 of the model is, the closer the model is to the true one.
  We propose then the following $k$-fold cross-validation procedure for model selection 
using the proposed estimate (\ref{estim div bis}) of $\varphi$-MI.
\begin{enumerate}
\item [(1)] Partition the sample  $(X_1,Y_1),\ldots,(X_n,Y_n)$ into $k$ equal size ($n_k$) subsamples. 
(Denote the $i$-th subsample $(X_{(i-1)n_k+1}, Y_{(i-1)n_k+1} ),   \ldots,(X_{in_k}, Y_{in_k})$, for all $i=1,\ldots,k$);
\item [(2)] Consider a candidate model $\cM_{\Theta_\ell}$; 
\item [(3)] From the sample $(X_1,Y_1),\ldots,(X_n,Y_n)$ remove the $i$-th subsample; compute the estimate $\widehat{\theta}_\ell^{(-i)}$ given by 
(\ref{estim param bis}) using the remaining $n-n_k$ observations, i.e.,
$$\widehat{\theta}_\ell^{(-i)} = \arg\sup_{\theta_\ell\in\Theta_\ell} M_{n-n_k}(\theta_\ell);$$ 
 \item [(4)] Repeat steps (2) and (3) for all $i=1,\ldots,k$, and obtain  
the following ``estimate'' 
$$C_V(\cM_{\Theta_\ell}):= \frac{1}{k}\sum_{i=1}^k \left(  \frac{1}{n_k} \sum_{j=(i-1)n_k+1}^{in_k}  f_{\widehat{\theta}_\ell^{(-i)}}(X_j,Y_j) -\frac{1}{n_k^2}
\sum_{j,m=(i-1)n_k+1}^{in_k} g_{\widehat{\theta}_\ell^{(-i)}}(X_j,Y_m) \right)   $$
of the expected criterion $M(\widehat{\theta}_\ell)$, i.e., 
$$M(\widehat{\theta}_\ell) = \int_{\Om} f_{\widehat{\theta}_\ell}(x,y)\, \der\P(x,y) - \int_{\Om} g_{\widehat{\theta}_\ell}(x,y)\, \der\P(x,y)^\perp;$$
\item [(5)]  Repeat steps (2-4) for all $\ell=1,\ldots,L$, and select the ``optimal'' model $\cM_{\Theta_{\ell^*}}$ that maximizes $C_V(\cM_{\Theta_\ell})$ 
over $\ell=1,\ldots, L$, i.e.,  
the model $\cM_{\Theta_{\ell^*}}$ with
$$\ell^* :=\arg\sup_{\ell\in\{1,\ldots,L\}} \, C_V(\cM_{\Theta_\ell}).$$
\end{enumerate}
Other model selection-type procedures can be investigated, through e.g.
correcting the bias of $M_n(\widehat{\theta}_\ell)$ as an estimate of the expected criterion 
$M(\widehat{\theta}_\ell)$, and selecting
the model that maximizes the obtained information criterion corrected from bias.
The correction  can be made e.g. by asymptotic evaluation of the bias
as in classical AIC criterion, or using bootstrap; see e.g.  \cite{Konishi_Kitagawa2008_Sringer} 
and \cite{Shao_Tu1995_Springer}.

\section{Asymptotic properties of the estimates} \label{SecAsympProp}
We state in Section \ref{Consistency} the consistency of both estimates $\widehat{I}_\varphi$ and $\widehat{\theta}_\varphi$, of the
$\varphi$-MI and the parameter $\theta_T$. Section \ref{limit distribution} gives, under the null hypothesis of independence, 
the limiting distribution
of the estimate $\widehat{I}_{\varphi_1}$ of the KL-MI,  as well as the corresponding 
estimate $\widehat{\theta}_{\varphi_1}$ of the parameter $\theta_T$, for some specific forms
of the model $\left\{h_\theta(\cdot,\cdot);\, \theta\in\Theta\right\}$.
Section \ref{Bootstrap calibration} provides bootstrap calibration 
of the critical value of any  $\widehat{I}_\varphi$-based test statistic
for general forms of the model $\left\{h_\theta(\cdot,\cdot);\, \theta\in\Theta\right\}.$ 

\subsection{Consistency} \label{Consistency}
In this section, we state consistency of the estimate $\widehat{I}_\varphi$, of the $\varphi$-MI, 
defined by (\ref{estim div}), as well as the consistency of the estimates 
$\widehat{\theta}_\varphi$ of $\theta_T$. 
We will use classical 
techniques from M-estimation theory. 
We will make use of the following conditions.
\begin{enumerate}
 \item[(A.5)] The parameter space $\Theta$ is a compact subset of $\R\times \R^d$ ;
 \item[(A.6)]  $\int_{\Om} \sup_{\theta\in\Theta}  \left|f_\theta(x,y)\right|\,\der\P(x,y)<\infty$;
 \item[(A.7)]  $\int_{\Om} \sup_{\theta\in\Theta}  g_\theta(x,y)^2\,\der\P^\perp(x,y)<\infty$,
\end{enumerate}
where $f_\theta$ and $g_\theta$ are defined respectively by~(\ref{la fonction f}) and~(\ref{la fonction g}).
Note that assumptions (A.6-7) imply  (A.3-4).

\begin{proposition} \label{MITPropConv1}
Assume that conditions (A.1, 5-7) hold. Then, the estimates $\widehat{I}_{\varphi}$ of $I_\varphi(\P)$ defined by~(\ref{estim div}) 
 and the estimates $\widehat{\theta}_\varphi$ of $\theta_T$ defined by~(\ref{estim param})  are consistent. 
Precisely, as $n\to\infty$, the following convergences in probability hold
$$\widehat{I}_\varphi \to I_\varphi(\P) \quad \text{and} \quad \widehat{\theta}_\varphi \to \theta_T.$$
\end{proposition}

\begin{remark}
Since in practice, all models are generally ``misspecified'', the true parameter value $\theta_T$ may not exist, it can however 
be replaced by the ``pseudo-true'' value 
$\theta_T^* := \arg\sup_{\theta\in\Theta} M(\theta),$
and the results of consistency in the above proposition remain valid.
\end{remark}

\subsection{The limiting distribution of the estimate $\widehat{I}_{\varphi_1}$ of KL-MI}\label{limit distribution}
We will give now the limiting distribution of the particular statistical test based on the estimate 
$\widehat{I}_{\varphi_1}$ of classical KL-MI, for specific forms of the model $h_\theta(\cdot,\cdot)$, under the null hypothesis of independence 
$\mathcal{H}_0 : \P=\P^\perp$.  
Consider the following specific form of the model $h_\theta(\cdot,\cdot)$ 
\begin{equation}\label{modele exp particulier}
h_\theta(x,y) = \exp\left(\alpha + m_\beta(x,y)\right)
\quad \text{with} \quad m_\beta(x,y) := \sum_{k=1}^d\beta_k \xi_k(x)\zeta_k(y),
\end{equation}
for some specified measurable real valued functions $\xi_k$ and $\zeta_k$, $k=1,\ldots,d$, defined, respectively, on $\mathcal{X}$ and $\mathcal{Y}$.
The parameter $\theta$ is the vector $\theta :=(\alpha,\beta_1,\ldots,\beta_d)^\top\in\Theta\subset\R\times\R^d.$
In this case, the functions (\ref{la fonction f}) and (\ref{la fonction g}) become
$$ f_\theta(x,y) = \alpha+\sum_{k=1}^d\beta_k \xi_k(x)\zeta_k(y)$$
and 
$$g_\theta(x,y) = \exp\left(\alpha +\sum_{k=1}^d\beta_k \xi_k(x) \zeta_k(y)\right)-1.$$ 
The value $\theta_0$, corresponding to the independence,  here is $\theta_0 = \boldsymbol{0}:=(0,\ldots,0)^\top\in \R^{1+d}$.
We will give the limiting distributions of $\widehat{\theta}_{\varphi_1}$ and $\widehat{I}_{\varphi_1}$,
 under the null hypothesis of independence $\P=\P^\perp$, i.e., when 
$\theta_T=\theta_0=\boldsymbol{0}.$  
We will consider the following assumptions.

\begin{enumerate}
 \item [(A.8)] There exists a neighborhood $N(\theta_T)$ of $\theta_T$ such that the third order partial derivative functions
 $\left\{ (x,y)\mapsto (\partial ^3/\partial^3\theta) f_\theta(x,y);\, \theta\in N(\theta_T) \right\}$ \\
 (resp. $\left\{ (x,y)\mapsto (\partial ^3/\partial^3\theta) 
 g_\theta(x,y);\, \theta\in N(\theta_T) \right\}$) are dominated by some functions $\P$-integrable (resp. some function $\P^\perp$-square-integrable);
 \item [(A.9)] The integrals  $\P \left\|f'_{\theta_T}\right\|^2$, $\P^\perp\left\|g'_{\theta_T}\right\|^2$, $\P \left\|f''_{\theta_T}\right\|$,
 $\P^\perp \left\|g''_{\theta_T}\right\|^2$ exist, and the matrix 
 \begin{equation} \label{MITEqnVar1}
  \Sigma_1:= -\left(\P f''_{\theta_T}-\P^\perp g''_{\theta_T} \right)
 \end{equation}
  is nonsingular.
\end{enumerate}

\begin{theorem} \label{MITThmConv2}
 Assume that conditions (A.1-2,5-9) hold and that $\P=\P^\perp$ (i.e., $\theta_T=\boldsymbol{0}$). Then, 
 \begin{enumerate}
  \item [(a)] $\sqrt{n}\, \widehat{\theta}_{\varphi_1}$ converges in distribution to a centered multivariate normal random variable 
  with covariance matrix  $\Sigma = \Sigma_1^{-1} \Sigma_2 \Sigma_1^{-1}$, where $\Sigma_1$ and $\Sigma_2$ are given respectively by~(\ref{MITEqnVar1}) and~(\ref{MITEqnVar2});

  \item [(b)] $2n\,  \widehat{I}_{\varphi_1}$ converges in distribution to the random variable $Z^\top Z$, where $Z$ is a centered multivariate
  normal random variable with covariance matrix 
  $$C = \Sigma_1^{-1/2} \Sigma_2 \Sigma_1^{-1/2}.$$
 \end{enumerate}
\end{theorem}

\begin{remark}
For the finite-discrete case, using the modeling (\ref{modele exponentiel fini bis}) in Example \ref{modele fini},
we can see that the corresponding matrix $\Sigma_2$ is of rank $(K_1-1)(K_2-1)$
and that the limiting distribution of $2n\widehat{I}_\varphi=2n\widehat{I}_\varphi^{emp}$  is a $\chi^2$-distribution 
with $(K_1-1)(K_2-1)$ degrees of freedom,
in particular, we recover the classical $\chi^2$-independence test theorem (for the case of finite-discrete distributions).
\end{remark}

\subsection{Bootstrap calibration}  \label{Bootstrap calibration}
In the general context of model (\ref{le modele}), for a given $\varphi$-MI, we propose the following
bootstrap procedure to calibrate the critical value of the corresponding test statistic. The critical value, denote it $b_\alpha$, 
 is the upper $\alpha$-quantile of the distribution  of the test statistic 
$S_n:=2n\widehat{I}_\varphi$, under the null hypothesis $\mathcal{H}_0$ of independence. 

\begin{enumerate}
\item Generate bootstrap sample $(X_1^*,Y_1^*),\ldots, (X_n^*,Y_n^*)$ from the product empirical distribution $\widehat{\P}^\perp=
\widehat{\P}_1\otimes \widehat{\P}_2$ of the original sample $(X_1,Y_1),\ldots, (X_n,Y_n)$;
\item Compute the value of the statistic $S_n^*:=2n\widehat{I}_\varphi^*$ from the bootstrap sample;
\item Repeat steps (1) and (2) $B=1000$ times, independently, to obtain the realizations $\left\{S_{n,1}^*,S_{n,2}^*,\ldots,S_{n,B}^*\right\}$;
\item Estimate $b_\alpha$ by $\widetilde{b}_\alpha:=$ the $(1-\alpha)$th quantile of the sequence 
$\left\{S_{n,1}^*,S_{n,2}^*,\ldots,S_{n,B}^*\right\}$. 
\end{enumerate}

\section{Large deviations principle and Bahadur asymptotic efficiency} \label{LDPBaha} \label{SecBahadurEff}
In this section, we compare Bahadur asymptotic efficiency of $\varphi$-MI based independence tests and show that the test based on classical 
Kullback-Leibler mutual information is the most efficient. 
Given $(\widehat{I}_{\varphi_1})_n$ and $(\widehat{I}_{\varphi_2})_n$ two sequences  of statistics,  for the test problem (\ref{le probleme de 
test}),  
numbers $\alpha \in (0,1)$, $\gamma \in (0,1)$ and an alternative hypothesis $\P\neq \P^\perp$, 
we define $n_i(\alpha, \gamma, \P)$, for $i \in \{1, 2\}$, respectively, as the 
minimal number of observations needed for the test based on ${\widehat{I}_{\varphi_i}}$ to have signification level $\alpha$ and power 
level $\gamma$. Then, Bahadur asymptotic relative efficiency of $(\widehat{I}_{\varphi_1})_n$ with respect to $(\widehat{I}_{\varphi_2})_n$ 
is defined 
as (if the limit exists)
$$\lim_{\alpha \rightarrow 0} \frac{n_2(\alpha, \gamma, \P)}{n_1(\alpha, \gamma, \P)}.$$
It is well known, see for example \cite{Nikitin1995_CUP} and   \cite{vanderVaart1998_CUP} Chapter 14, 
that if both sequences $(\widehat{I}_{\varphi_1})_n$ and $(\widehat{I}
_{\varphi_2})_n$ satisfy a large deviation principle under the null hypothesis (with good rate functions $e_{\varphi_1}(\cdot)$ and $e_{\varphi_2}(\cdot)
$) 
and also a law of large number under a given alternative hypothesis $\mathcal{H}_1 : \P\neq \P^\perp$, 
with asymptotic means $\mu_{\varphi_1}(\P)$ and $\mu_{\varphi_2}(\P)$, respectively, 
then the Bahadur asymptotic relative efficiency equals $e_{\varphi_1}(\mu_{\varphi_1}(\P)) / e_{\varphi_2}(\mu_{\varphi_2}(\P))$. Particularly, the most 
efficient test maximizes  Bahadur slope $e_{\varphi}(\mu_{\varphi}(\P))$.
A law of large number under the alternative hypothesis is given for the sequence $(\widehat{I}_{\varphi})_n$ in Proposition
\ref{MITPropConv1} above; the excepted value $\mu_{\varphi}(\P)$ being $\mu_{\varphi}(\P) = I_{\varphi}(\P) =D_{\varphi}(\P,\P^\perp)$. 
The following theorem establishes a large deviation principle under the null hypothesis of independence. 
It relies on some generalization due to \cite{Eischelsbacher_Scmock2002_AIHP} 
of classical Sanov theorem to finer topologies and the contraction principle.
Let $\mathcal{G}$ be the set of measurable functions, from $\Om$ into $\mathbb{R}$, 
given by
$$\mathcal{G} := \mathcal{B} \cup \{ \varphi'(h_{\theta}); \theta \in \Theta\} \cup \{ \varphi^*(\varphi'(h_{\theta})); \theta \in \Theta\},$$
where $\mathcal{B}$ is the set of all measurable bounded functions from $\Om$ into $\R$.
Recall that $\mathcal{M}_1 = \mathcal{M}_1(\Om)$ is the set of all probability measures on $\Om$, and let us
 introduce the subset $$\mathcal{M}_\mathcal{G}:=\mathcal{M}_\mathcal{G}(\Om) 
 := \left\{P \in \mathcal{M}_1:   \int_{\Om} |\varphi'(h_\theta)| \, dP < \infty, \int_{\Om} |\varphi^*(\varphi'(h_\theta))| \, dP^{\perp}< \infty, 
 \forall \theta \in \Theta 
 \right\}.$$ 
Define on $\mathcal{M}_\mathcal{G}$ the $\tau_\mathcal{G}$-topology as the coarsest one that makes applications 
$P \in \mathcal{M}_\mathcal{G} \mapsto \int_{\Om} \varphi'(h_\theta) \, dP$, $P \in \mathcal{M}_\mathcal{G} \mapsto \int_{\Om} f \, dP$, 
$P \in \mathcal{M}_\mathcal{G} \mapsto \int_{\Om} \varphi^*(\varphi'(h_\theta)) \, dP^\perp$ and $P \in \mathcal{M}_\mathcal{G} \mapsto 
\int_{\Om} f \, dP^\perp$   continuous, 
for all $\theta \in \Theta$ and all $f \in \mathcal{B}$.
Finally, define, for all  $Q \in \mathcal{M}_\mathcal{G}$, the ``pseudo-divergence''
$$\Ddp(Q,Q^{\perp}):=\sup_{\theta \in \Theta} \left\{ \int_{\Om} \varphi'(h_{\theta}(x,y))\, dQ(x,y) - 
\int_{\Om} \varphi^*\left( \varphi'(h_{\theta}(x,y))\right) 
dQ^{\perp}(x,y) \right\}.$$
Obviously, $\Ddp(Q,Q^{\perp}) \leq \Dp(Q,Q^{\perp}) =: I_{\varphi}(Q)$ with equality for probability distributions such that $dQ/ dQ^{\perp} = 
h_{\theta}$ for some $\theta \in \Theta$. Note also that $Q\in \mathcal{M}_\mathcal{G}\mapsto \Ddp(Q,Q^{\perp})$
 is continuous with respect to the $\tau_\mathcal{G}$-topology as the 
supremum over the compact set $\Theta$ 
of continuous functions.
The large deviation principle for the sequence $(\widehat{\P}(\cdot))_{n}$ of empirical measures defined by (\ref{MITEqnProbEmp}), 
established by
\cite{Eischelsbacher_Scmock2002_AIHP}, requires the existence of exponential moments; in the context of the model~(\ref{le modele}), we 
thus 
assume
\begin{itemize}
 \item[(A.10)] for all $f \in \mathcal{G}$, for all $a>0$, 
 $$\quad\int_{\Om} \exp(a |f|) \,  \der \P< \infty.$$
\end{itemize} 

Note that the strong assumption $(A.10)$ implies (A.3-4) if $\P=\P^\perp$.
In the context of the models described in Examples~\ref{modele gaussien} to~\ref{modele copule nonparam}, assumption (A.10) may not be satisfied for 
some $\varphi$-divergences ; particularly, it does not generally hold for power-divergences (except for finite-discrete distribution models described in 
Example \ref{modele fini}). 
A sufficient condition for (A.10) is 
\begin{itemize}
 \item[(A.11)] there exist real numbers $m,M \in (a_{\varphi}, b_{\varphi})$ such that
 $m<h_{\theta}(x,y) <M$, \, $\forall(x,y) \in \Om,\,  \forall \theta\in\Theta$.
\end{itemize}
Indeed, for all $a>0$, the functions $\exp(a|\varphi'(h_{\theta})|)$ and $\exp(a|\varphi^*(\varphi'(h_{\theta}))|)$ 
are bounded and therefore integrable with respect to both
$\P$ and $\P^{\perp}$.
Again, (A.11) is not generally satisfied for models described in the previous examples for power-divergences, but it may 
be artificially verified by truncating the distributions in the models. 
Let us also point out that Theorems~\ref{LDP} and~\ref{Baha} below may remain true with some alternative assumptions on the distribution queues, lighter than 
(A.10). Particularly, simulations performed in Section~\ref{SecSimus} for bivariate Gaussian distributions tend to show that Theorem \ref{Baha} holds for the 
Gaussian model described in Example~\ref{modele gaussien}. 
For getting a closed form for the LDP of $(\widehat{I}_\varphi)_n$, we will establish the right-continuity of the rate function, making use of one of the 
following assumptions:
\begin{enumerate}
 \item[(A.12.a)] $(X,Y)$ is finite-discrete, supported by $\Om$;
 \item[(A.12.b)] The model $\left\{h_\theta(\cdot,\cdot); \, \theta =(\alpha,\beta^\top)^\top\in \Theta\right\}$ is of the from 
 $h_\theta(x,y)=\exp\left(\alpha+m_\beta(x,y)\right)$ with the condition that,  for any constant $c$ and  any $\beta$,  we have
  $\P^\perp \left(m_\beta(X,Y) = c\right) \neq 0$ iff  $\beta = (0,\ldots,0)^\top$ and $c=0$.
\end{enumerate}

\begin{theorem} \label{LDP}
 Let $(X,Y)$ be a couple of independent random variables with joint distribution $\P=\P^{\perp} \in \cM_{\Theta} \cap \mathcal{M}_\mathcal{G}$. 
 \begin{enumerate}[(1)]
\item Suppose that 
 conditions (A.1-2, 5-7, 10 and 12.b) are satisfied.
 Then, the sequence $(\widehat{I}_{\varphi})_n$ of estimates, of $I_{\varphi}(\P)=0$, given by (\ref{estim div}), satisfies the following large deviation principle 
 \begin{equation} \label{BEITEqnLDP}
  \frac{1}{n} \log \P^\perp \left( \widehat{I}_{\varphi} >d \right) \stackrel{n \rightarrow \infty}{\longrightarrow} - e_{\varphi}(d), \quad d>0,
 \end{equation}
 where the good rate function $e_{\varphi}(\cdot)$ is 
 \begin{equation} \label{BEITeqnFTLDP}
  e_{\varphi}(d):=\inf_{Q\in\Omega_d} \mathbb{K}(Q, \P^{\perp})  \, \text{ with } \,  \Omega_d := \left\{ Q\in \mathcal{M}_\mathcal{G} \text{ such that } 
  \Ddp(Q,Q^{\perp}) \geq d\right\}.
 \end{equation}
 \item Assume that conditions (A.1-2, 5 and 12.a) are satisfied. Then the above statement holds if $\mathcal{M}_\mathcal{G}$
 is replaced by the set of all discrete-finite distributions with the same finite support $\Om$.
 \end{enumerate}
\end{theorem}

\vspace{0.2cm}
In view of Proposition \ref{MITPropConv1} and Theorem \ref{LDP} above, 
the Bahadur slope of the independence test based on $\widehat{I}_{\varphi}$, for any $\varphi$, is given then by
\begin{eqnarray*}
 s_{\varphi} &: = & e_{\varphi}(I_{\varphi}(\P) )\\
  &=& \inf\{ \mathbb{K}(Q, \P^{\perp}): \Ddp(Q,Q^{\perp}) \geq \Dp(\P,\P^{\perp})\}.
\end{eqnarray*}
Since $\Ddp(\P,\P^{\perp}) = \Dp(\P,\P^{\perp})$, we have $\P \in \{Q: \Ddp(Q,Q^{\perp}) \geq D_{\varphi}(\P,\P^{\perp})\}$, so that, for any $\varphi$,
\begin{equation} \label{BEITinqBSphi}
 s_{\varphi} \leq \mathbb{K}(\P , \P^{\perp}) = I_{KL}(\P) = I_{\varphi_1}(\P).
\end{equation}
Equality is achieved in (\ref{BEITinqBSphi}) for the divergence $D_{\varphi} = \K$. Indeed,
$$s_{KL} = \inf\{ \mathbb{K}(Q,\P^{\perp}): \mathcal{D}_{KL}(Q,Q^{\perp}) \geq \mathbb{K}(\P,\P^{\perp})\}.$$
Straightforward computations yield 
$$\mathbb{K}(Q,\P^{\perp}) = \mathbb{K}(Q,Q^{\perp}) + \mathbb{K}(Q_1,\P_1) + \mathbb{K}(Q_2,\P_2),$$
for any $Q \in \mathcal{M}_\mathcal{G}$. 
Particularly, for any $Q$ such that $\mathcal{D}_{KL}(Q,Q^{\perp}) \geq \mathbb{K}(\P,\P^{\perp})$, we have $\K(Q, Q^\perp) \geq 
\mathcal{D}_{KL}(Q,Q^{\perp}) \geq \mathbb{K}(\P,\P^{\perp})$, hence, 
$$\mathbb{K}(Q,\P^{\perp}) \geq \mathbb{K}(Q,Q^{\perp}) \geq \mathbb{K}(\P,\P^{\perp}),$$
so that 
\begin{equation}\label{BEITinqBS KL}
s_{KL} \geq \mathbb{K}(\P,\P^{\perp}).
\end{equation}
Combining (\ref{BEITinqBSphi}) and (\ref{BEITinqBS KL}), we obtain 
\begin{theorem} \label{Baha}
 Let $(X,Y)$ be a couple of random variables with joint distribution $\P \in \mathcal{M}_{\Theta} \cap 
 \mathcal{M}_\mathcal{G}$.  
 Suppose that either conditions (A.1-2, 5-7, 10 and 12.b) or (A.1-2, 5 and 12.a) are satisfied.
 For the test problem (\ref{le probleme de test}), the test based on the estimate $\widehat{I}_{\varphi_1}$, see 
 (\ref{estim div}), of the Kullback-Leibler mutual information, is uniformly (i.e., whatever be the alternative 
 $\P \neq \P^\perp$) the most efficient test, in Bahadur sense, among all $\widehat{I}_{\varphi}$-based tests, 
 including the classical $\chi^2$-independence one.
\end{theorem}

\begin{remark}
Assume that $\P$ is a finite-discrete distribution.
We obtain then that KL-MI based independence test is more efficient than 
the classical $\chi^2$ independence one. This result was already stated, in goodness-of-fit testing 
for finite-discrete distributions, 
see e.g.  \cite{vanderVaart1998_CUP} Chapter 17 Section 17.6. The above 
theorem extends it to testing independence, for more general probability distributions, 
not necessarily finite-discrete.
\end{remark}

\section{Simulations} \label{SecSimus}

This Section aims at numerically comparing through simulations $\varphi$-MI based tests with other independence or non-correlation tests. Precisely, 
Section~\ref{SecSimusFinite} focuses on finite-discrete random vectors, for which the optimal KL-MI test is compared to the very popular 
(but not optimal) 
$\chi^2$-independence test. Section~\ref{SecSimusGauss} compares KL-MI and $\chi^2$ tests to classical non-correlation tests of Pearson, Kendall 
and Spearman. Finally, Section~\ref{SecSimusCopula} deals with the example of the copula density model of Farlie-Gumbel-Morgenstern (FGM), for which 
the critical values of KL-MI and $\chi^2$-MI tests are derived through the bootstrap procedure described in Section \ref{Bootstrap calibration}.

\subsection{Testing independence of finite-discrete random variables}\label{SecSimusFinite}
\begin{figure}[ht!]
 \centering
  \includegraphics[width=0.7\textwidth]{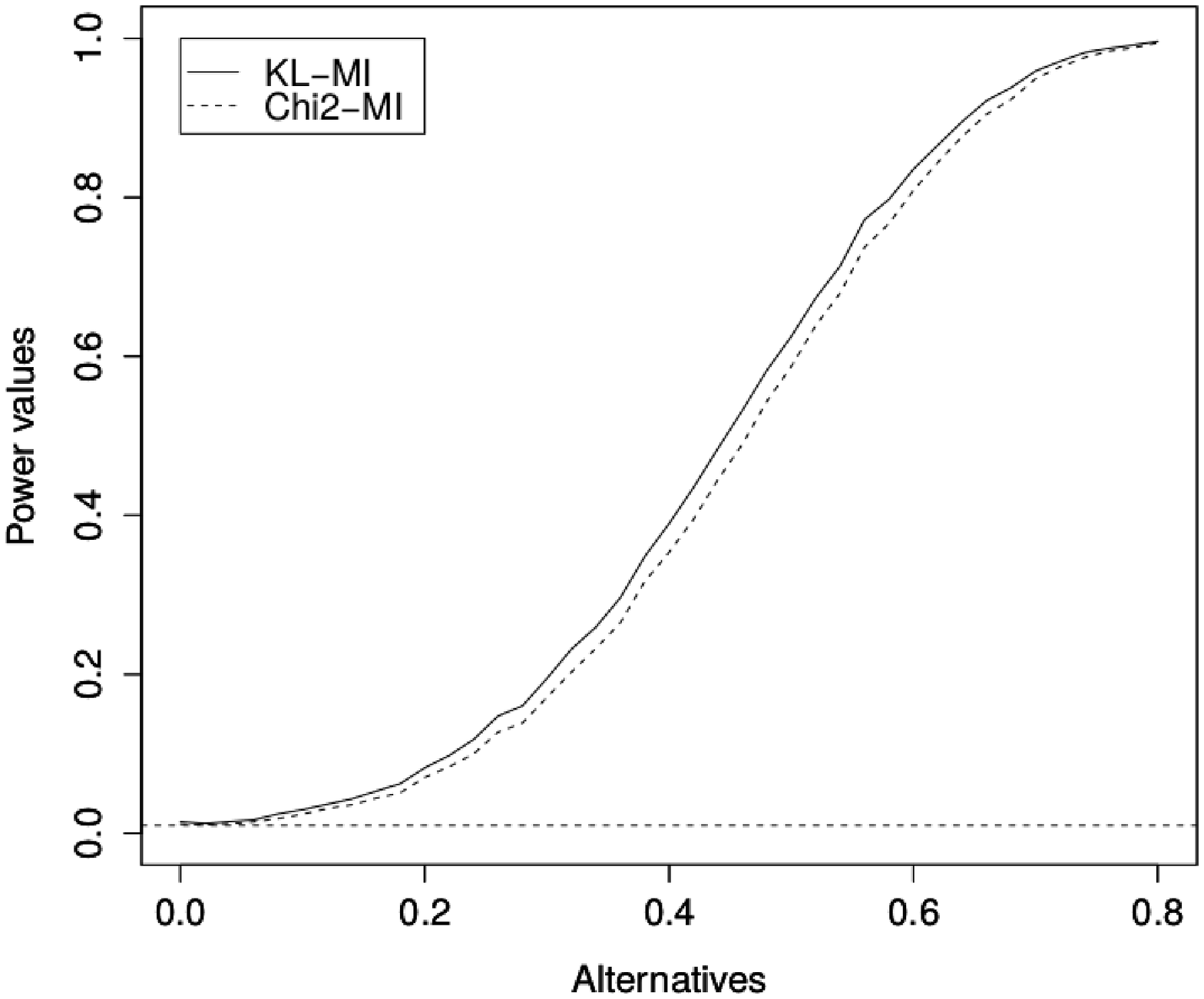}
 \caption[]{Comparison of KL-MI and $\chi^2$-MI based tests for finite-discrete random variables taking values in $\{1,2\}$, with $n=30$.}
 \label{figure_KL_Chi2_2classes}
\end{figure}

\begin{figure}[ht!]
 \centering
  \includegraphics[width=0.7\textwidth]{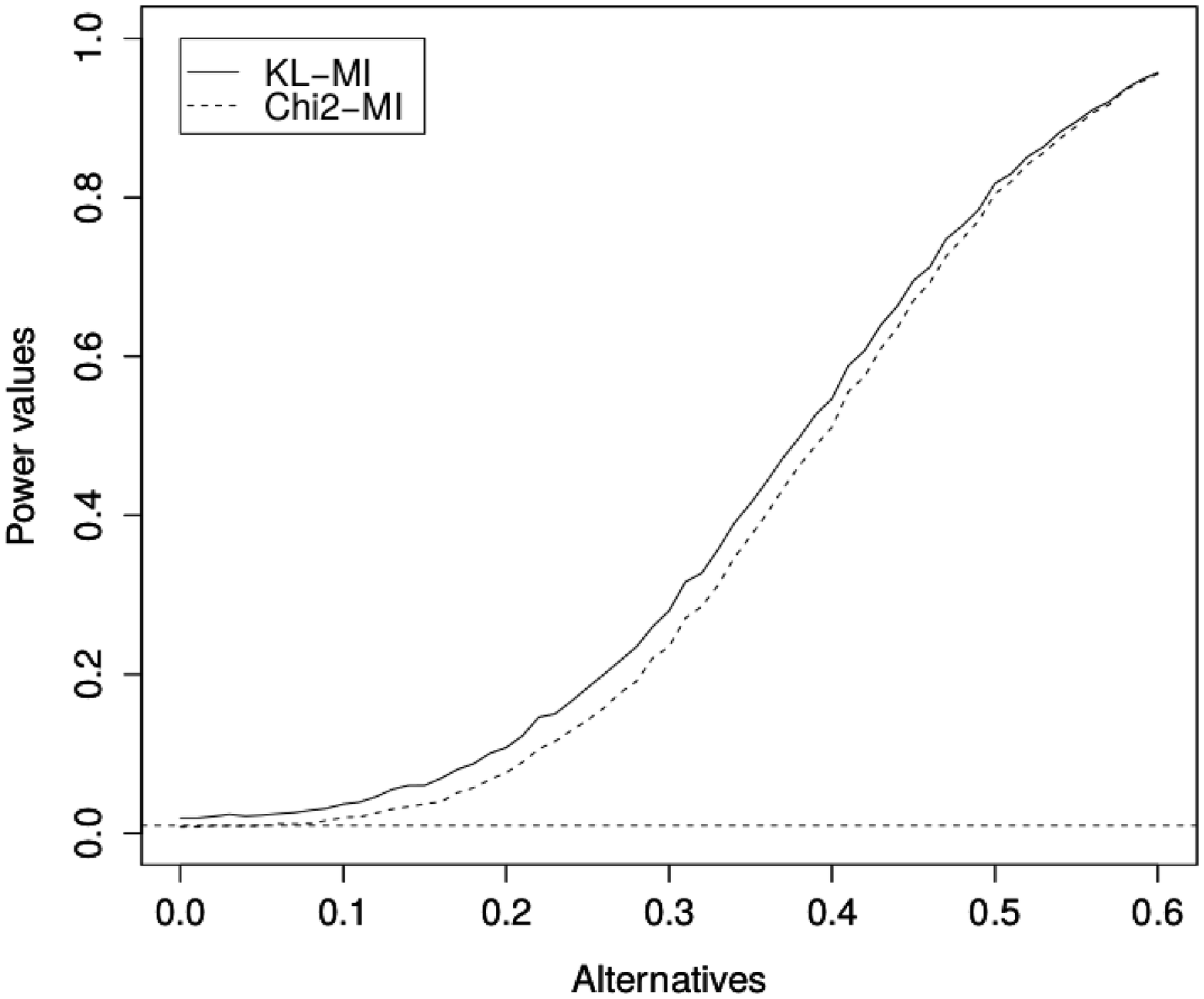}
 \caption[]{Comparison of KL-MI and $\chi^2$-MI based tests for finite-discrete random variables taking values in $\{1,2,3\}$, with $n=35$.}
 \label{figure_KL_Chi2_3classes}
\end{figure}

As stated in Example~\ref{MITExpDualEstFinite}, the dual estimates $\widehat{I}_{\varphi}$ given by (\ref{estim div}) equal the direct empirical ones 
(\ref{estimMIempDirect}). 
Their properties and asymptotic behavior are well-known; see e.g. \cite{Pardo2006_CH}. They are recovered by Propositions \ref{MITPropConv1}, 
Theorem \ref{MITThmConv2} and Theorem~\ref{Baha}.
We illustrate these properties through simulations, by comparing the power of KL-MI and $\chi^2$-MI tests, 
for various sample sizes and finite-discrete supports  
$\mathcal{X}=\mathcal{Y}=\left\{1,\ldots,K\right\}$, 
and for alternatives $P\in \cM_1(\Om)$ of the form $P_\theta :=(p_{x,y;\theta})_{(x,y)}$, with
\begin{equation} \label{MITEqnDistFiniteSimus}
 p_{x,y;\theta} =(1-\theta) \frac{1}{K^2}+ \theta \frac{1}{K} \mathds{1}_{\{x=y\}}, \quad (x,y) \in \Om,
\end{equation}
where $K=|\mathcal{X}| = |\mathcal{Y}|$ and $\theta \in (0,1)$, i.e., the random variables $X$ and $Y$ are uniformly distributed on the set
$\left\{1,\ldots,K\right\}$,  and the conditional 
distribution $P_{Y|X=x}(\cdot)$, of $Y$ knowing $X=x$, is the mixture of the uniform distribution on $\left\{1,\ldots,K\right\}$ with weight $(1-\theta)$
and the Dirac measure $\delta_x(\cdot)$ with weight $\theta$, for all $x\in\left\{1,\ldots,K\right\}.$
 Hence, for $\theta = \theta_0=0$, $X$ and $Y$ are independent, while for $\theta=1$, we have $Y =X$.
The level of the tests has been set to $\alpha=0.01$. 
The asymptotic distribution of $2n\widehat{I}_\varphi$ is $\chi^2\left((K-1)(K-1)\right)$, a $\chi^2$-distribution with $(K-1)^2$ degrees of freedom,
for both KL-MI or $\chi^2$-MI. The critical value  $b_{0.01}$ of both test statistics  is taken then to be the upper $0.01$-quantile 
of the $\chi^2\left((K-1)(K-1)\right)$-distribution.
 Then, we have estimated their respective 
powers, by means of Monte-Carlo procedure from $10000$ samples drawn according to $P_\theta$ 
given by (\ref{MITEqnDistFiniteSimus}), 
for various  mixture parameter values $\theta \in (0,1)$. The results are presented in Table \ref{MITTablePowerFinite}, Figure 
\ref{figure_KL_Chi2_2classes} and Figure \ref{figure_KL_Chi2_3classes}. 
We can see that the KL-MI test outperforms the classical $\chi^2$ one. The nominal levels of both KL-MI and $\chi^2$-MI test statistics
are both close to the test level $\alpha=0.01$. 
\begin{table}[ht!]
 \centering
 \tiny
 \begin{tabular}{|lr|c|c|c|c|c|c|c|c|}
  \hline
  $K=|\mathcal{X}|=|\mathcal{Y}|=2$ & $\theta=$ & 0 & 0.08 & 0.18 & 0.28 & 0.38 & 0.48 & 0.58 & 0.68 \\ \hline
  \multirow{2}{*}{$n=30$} & KL-MI test power & 0.0123 & 0.0242 & 0.0647 & 0.1681 & 0.3343 & 0.5690 & 0.7981 & 0.9415 \\ 
   & $\chi^2$ test power  & 0.0102 & 0.0200 & 0.0550 & 0.1433 & 0.2968 & 0.5330 & 0.7703 & 0.9288 \\ \hline
  \multirow{2}{*}{$n=40$} & KL-MI test power & 0.0119 & 0.0213 & 0.0764 & 0.2176 & 0.4502 & 0.7180 & 0.9046 & 0.9850 \\ 
   & $\chi^2$ test power  & 0.0100 & 0.0184 & 0.0694 & 0.2006 & 0.4272 & 0.6970 & 0.8957 & 0.9839 \\ \hline \hline
   
  $K=|\mathcal{X}|=|\mathcal{Y}|=3$ & $\theta=$ & 0 & 0.07 & 0.15 & 0.23 & 0.31 & 0.39 & 0.47 & 0.55 \\ \hline
  \multirow{2}{*}{$n=35$} & KL-MI test power& 0.0192 & 0.0261 & 0.0604 & 0.1503 & 0.3162 & 0.5267 & 0.7476 & 0.8952 \\ 
   & $\chi^2$ test power  & 0.0081 & 0.0118 & 0.0371 & 0.1157 & 0.2708 & 0.4878 & 0.7259 & 0.8895 \\ \hline
  \multirow{2}{*}{$n=50$} & KL-MI test power & 0.0152 & 0.0261 & 0.0782 & 0.2152 & 0.4369 & 0.7150 & 0.9039 & 0.9816 \\ 
   & $\chi^2$ test power  &  0.0088 & 0.0167 &  0.0648  & 0.1929 & 0.4283 & 0.7124 & 0.9057 & 0.9832 \\ \hline 
 \end{tabular}
 \vskip 0.2cm
 \caption{Comparison of powers of KL-MI and $\chi^2$-MI tests. 
 The number of cells $K$ is indicated at the top left of each block.
The sample sizes $n$ are given by the first column 
 while the mixture parameter values $\theta$, see its definition in (\ref{MITEqnDistFiniteSimus}), are given by the first row.}
 \label{MITTablePowerFinite}
\end{table}

\subsection{Comparison of $\varphi$-MI based and noncorrelation tests in the Gaussian setting} \label{SecSimusGauss}
\begin{figure}[ht!]
 \centering
  \includegraphics[width=0.6\textwidth, angle=270]{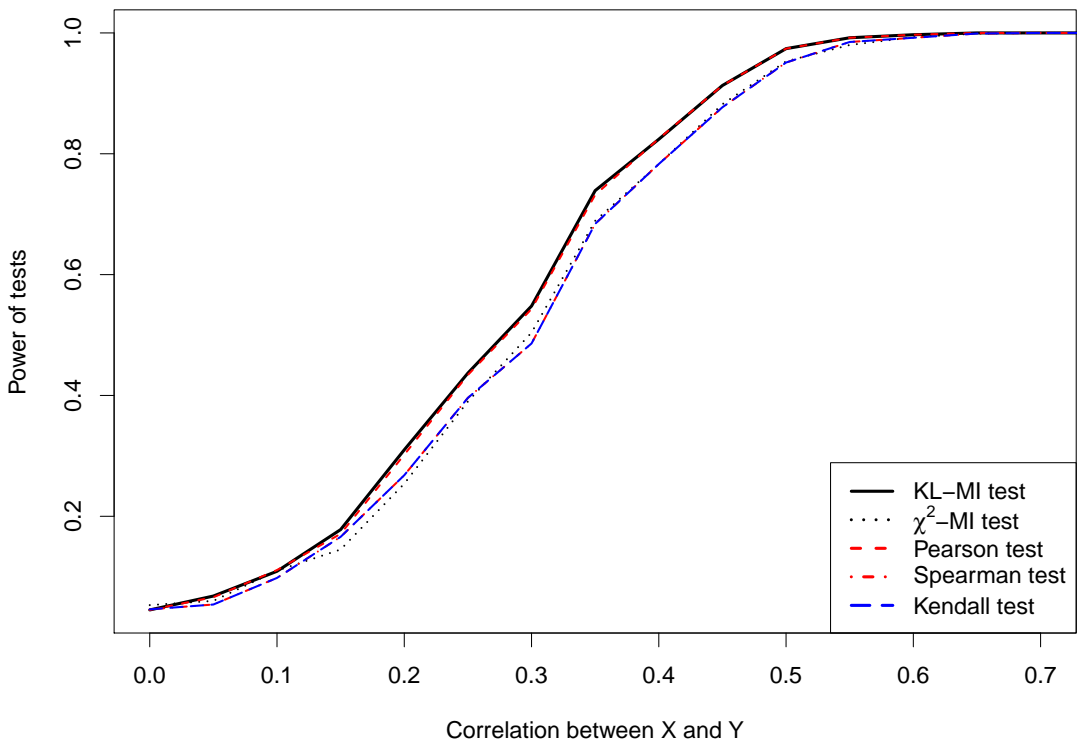}
 \caption[]{Comparison of powers of KL-MI and $\chi^2$-MI tests with noncorrelation tests of Pearson, Spearman and Kendall.}
 \label{MIITFigPower}
\end{figure}

For bidimensional normally distributed random vectors, the corresponding model $h_\theta(\cdot,\cdot)$, see Example \ref{modele gaussien}, 
is of the form (\ref{modele exp particulier}),
 so that the asymptotic 
distribution of the dual  KL-MI based test statistic $2n\widehat{I}_{\varphi_1}$ 
 is explicit. Hence, explicit (asymptotic) critical value can be obtained for the test statistic $2n\widehat{I}_{\varphi_1}$. Although assumption (A.10) may not 
be satisfied without restricting the support of $(X,Y)$ to a bounded subset of $\R^2$, we can compare numerically the powers of the $\varphi$-MI 
based tests. 
Precisely, in this Section we manage to  compare the powers of KL-MI and $\chi^2$-MI independence tests with noncorrelation tests for samples of 
size 
$n=50$ drawn according to bivariate normal distributions. We have fixed the level $\alpha=0.05$ and computed the critical value of 
KL-MI  based test 
by means of Monte-Carlo simulations of the asymptotic distribution of $2n\widehat{I}_{\varphi_1}$ given by Theorem \ref{MITThmConv2} 
(10000 samples of the variable $Z$ in Theorem~\ref{MITThmConv2} have been simulated; the critical value has been obtained as the $0.95$-quantile 
of the linearly interpolated empirical cumulative density function). The critical value for the $\chi^2$-MI based test have been estimated directly by 
simulating 10000 samples of size 50 of a bivariate Gaussian random vector with independent centered and reduced distribution and computing the 
$0.95$-quantile of the corresponding tes statistic $2n\widehat{I}_{\varphi_2}$.
Then we have estimated the power of these tests as well as noncorrelation tests of Pearson, Spearman and Kendall, 
still by Monte-Carlo methods: for any correlation value $\rho \in \{0,1/20, 
2/20, \dots, 1\}$, we have considered $N=1000$ samples, with size $n=50$, 
of centered bivariate Gaussian couples with marginal variances equal to $1$ and covariance $\rho$ 
varying from $0$ to $1$.
Recall that the noncorrelation test of Pearson, for this particular Gaussian model,  
is the most uniformly powerful test, among all tests with the same level $\alpha$.
Figure \ref{MIITFigPower} presents the power curves for KL-MI  (plain black curve), 
$\chi^2$-MI (dotted black curve) independence tests, and Pearson 
(dashed red curve), Kendall and Spearman (mixed dashed and dotted red and blue curves) 
correlation tests, obtained from $N=1000$ samples of size 
$n=50$ of bivariate Gaussian distributions. 
For this setting, we can see form Figure \ref{MIITFigPower}, that
our poposed KL-MI independence test is almost as powerful as the most uniformly powerful 
independence test of Pearson. $\chi^2$-MI, Spearman and Kendall tests 
have comparable powers, lower than KL-MI and Pearson's ones.

\subsection{Comparison of $\varphi$-MI based tests for a copula density model} \label{SecSimusCopula}
This Section aims at comparing numerically the $\varphi$-MI based independence tests in the context of 
semiparametric copula-type model, as described 
in Example \ref{modele copule nonparam}. 
We consider here the Farlie-Gumbel-Morgenstern (FGM) copula model
$$C_{FGM}(u,v;\theta) = uv(1+\theta (1-u)(1-v)), \quad (u,v) \in [0,1]^2, \quad \theta \in \Theta = [-1,1],$$
with $\theta_0=0$.
We compare the powers of KL-MI and $\chi^2$-MI based tests of independence to noncorrelation ones.
We consider the alternative hypothesis that $X$ and $Y$ are uniformly distributed on $[0,1]$ and copulated by a 
FGM copula. We consider values of the parameter $\theta$ of the form $\theta = k/16$, with $k \in \{0, \dots, 16\}$. 
We have estimated the critical values of the KL-MI and $\chi^2$-MI tests using the bootstrap procedure presented in Section \ref{Bootstrap calibration}, from an original sample of size $n=50$ resampled 10 000 times.
The powers are computed by Monte-Carlo method from 
$N=5000$ samples of size $n=50$. The results  are presented in Table~\ref{MITTablePowerFGM}.
We can see again that KL-MI based test still outperforms the others.  We can see also that the nominal levels 
(of KL-MI and $\chi^2$-MI test statistics) are sufficiently close to the 
test levels evaluated through the bootstrap procedure described in Section \ref{Bootstrap calibration},  
with  $\alpha=0.05$. 

\begin{table}[!ht]
\tiny
\centering
\begin{tabular}{|l|c|c|c|c|c|c|c|c|}
  \hline 
 $\theta$ & 0 & 1/16 & 2/16 & 3/16 & 4/16 & 5/16 & 6/16 & 7/16 \\ 
  \hline \hline
     KL-MI & 0.062 & 0.061 & 0.064 & 0.076 & 0.093 & 0.120 & 0.142 & 0.171\\ \hline
  $\chi^2$ & 0.054 & 0.055 & 0.057 & 0.066 & 0.084 & 0.108 & 0.129 & 0.160\\ \hline
   Pearson & 0.052 & 0.057 & 0.061 & 0.072 & 0.089 & 0.113 & 0.135 & 0.170\\ \hline
  Spearman & 0.055 & 0.058 & 0.060 & 0.069 & 0.086 & 0.110 & 0.133 & 0.164\\ \hline
   Kendall & 0.056 & 0.057 & 0.057 & 0.069 & 0.086 & 0.111 & 0.130 & 0.161\\ \hline
\end{tabular}

\vspace{0.2cm}
\begin{tabular}{|l|c|c|c|c|c|c|c|c|c|}
  \hline 
 $\theta$ & 8/16 & 9/16 & 10/16 & 11/16 & 12/16 & 13/16 & 14/16 & 15/16 & 1 \\ 
  \hline \hline
     KL-MI  & 0.219 & 0.261 & 0.312 & 0.382 & 0.431 & 0.498 & 0.565 & 0.622 & 0.691 \\ \hline
  $\chi^2$  & 0.202 & 0.244 & 0.296 & 0.362 & 0.404 & 0.472 & 0.527 & 0.589 & 0.659 \\ \hline
   Pearson  & 0.213 & 0.257 & 0.309 & 0.375 & 0.427 & 0.493 & 0.549 & 0.611 & 0.677 \\ \hline
  Spearman  & 0.207 & 0.249 & 0.300 & 0.369 & 0.410 & 0.478 & 0.533 & 0.596 & 0.663 \\ \hline
   Kendall  & 0.203 & 0.243 & 0.293 & 0.356 & 0.405 & 0.467 & 0.527 & 0.584 & 0.647 \\ 
   \hline
\end{tabular}
 \vspace{0.2cm}
 \caption{Power functions of KL-MI and $\chi^2$-MI tests 
 compared to noncorrelation tests obtained from $N=5000$ samples of size $n=50$ of the FGM 
 copula with parameter $\theta$ varying from 0 to 1 by step of $1/16$.} 
 \label{MITTablePowerFGM}
\end{table}

\section{Concluding remarks and discussion}

In this paper, we have defined and studied estimates of $\varphi$-mutual informations, based on the dual representation of 
$\varphi$-divergences and a 
semiparametric modeling of the density ratio between the joint distribution of the couple and the product distribution of its margins. The 
consistency of these estimates -- named dual-estimates -- has been established assuming some classical regularity conditions on the model; 
the asymptotic normality has been established for classical Kullback-Leibler mutual information  and specific models by means of classical 
M-estimation theory arguments. The asymptotic normality of other $\varphi$-mutual information dual-estimates may be derived similarly, for 
specific models depending on the considered $\varphi$-divergence. For example, when dealing with the power divergence associated to 
$\varphi_{\gamma}$ functions given by (\ref{MITEqnPowerFunction}), 
the asymptotic normality of the corresponding $\varphi_{\gamma}$-mutual-information 
dual-estimates may be derived in a similar way when focusing on the so-called $\gamma$-exponential semiparametric model
\begin{equation*}
\P \in \left\{P \in \cM_1(\Om) \text{ such that } 
 \frac{dP}{dP^{\perp}}(x,y) = 
 \exp_{\gamma} \left( \sum_{k=0}^d \theta_k \xi_k(x) \zeta_k(y) \right),  \theta=(\theta_0, \dots, \theta_d) \in \Theta \right\},
 \end{equation*}
where $\exp_\gamma(t) := ((\gamma-1)t+1)_+^{\frac{1}{\gamma-1}}$, with $(\cdot)_+=\max(0,\cdot)$.
Our semiparametric approach for estimating mutual informations constitutes a promising alternative to classical nonparametric procedures 
based on kernel density estimation or adaptive partitioning. No parameters such as bandwidth or kernel type has to be adjusted. The 
asymptotic normality of dual-estimates is also of significative  importance, particularly, for hypothesis-testing purpose. 
For the sake of both completeness and accessibility, we are developing a package for the R software providing user-ready procedures, 
including the $k$-fold cross validation procedure  described in Section~\ref{Sec model select}, for selecting the  model that best matches the 
data. 
We also aim at comparing the dual-estimates of mutual informations to nonparametric estimates.
As an application of dual-estimation of mutual informations, we have derived a class of independence tests, recovering as a particular case, 
the classical $\chi^2$-independence test. For a large variety of situations including finite-discrete random couples, the most efficient test is 
based on the KL-MI estimates, outperforming the classical $\chi^2$-independence one. Motivated by the simulation experiments presented 
in this paper, we guess that the optimality of KL-MI independence test can be extended to a larger family of models.

\section{Appendix}
\textit{Proof of Proposition \ref{MITPropConv1}}. 
Using continuity of $g_\theta(x,y)$ in $\theta$ on the compact set $\Theta$, and condition (A.7), we can 
state, by Bienaym\'e-Tchebychev 
inequality, the 
uniform convergence in probability 
\begin{equation}\label{cv 1}
A_n:=\sup_{\theta\in\Theta} \left|  \int g_\theta(x,y)\, \der\widehat{\P}^\perp(x,y) - \int g_\theta(x,y)\, \der\P^\perp(x,y) \right| \to 0.
\end{equation}
Under condition (A.6), using continuity of $f_\theta(x,y)$ in $\theta$ over the compact set $\Theta$, we have by uniform weak law of large numbers 
the convergence in probability 
\begin{equation}\label{cv 2}
B_n:=\sup_{\theta\in\Theta} \left|  \int f_\theta(x,y)\, \der\widehat{\P}(x,y) - \int f_\theta(x,y)\, \der\P(x,y) \right| \to 0.
\end{equation}
 Now, we have 
 \begin{eqnarray}
  \left| \widehat{I}_\varphi - I_\varphi(\P) \right| & =  & \left| \sup_{\theta\in\Theta} M_n(\theta)-\sup_{\theta\in\Theta} M(\theta)   \right| \nonumber\\
  & =  &  \left| M_n(\widehat{\theta}_\varphi)-M(\theta_T)   \right|  : = 
  \left|C_n\right| 
 \end{eqnarray}
  with $$ C_{n,L} := M_n(\theta_T) -M(\theta_T)  \leq C_n \leq M_n(\widehat{\theta}_\varphi)-M(\widehat{\theta}_\varphi) =: C_{n,R}.$$
We can see that both sides converge in probability to zero, since 
$$\left|  C_{n,L} \right| \leq A_n+B_n  \quad \text{and} \quad \left|  C_{n,R} \right| \leq A_n+B_n $$
and the use of convergences (\ref{cv 1}) and (\ref{cv 2}). 
We conclude that  $\widehat{I}_\varphi \to I_\varphi(\P)$ in probability.
The convergence of $\widehat{\theta}_\varphi$ to $\theta_T$ holds by direct application of Theorem 5.7 in \cite{vanderVaart1998_CUP},
using the uniform convergence in probability
$$\sup_{\theta\in\Theta} \left|M_n(\theta)-M(\theta)\right|\to 0$$
and the well-separability of the supremum $\theta_T$; it is unique and interior point of 
$\Theta$.  \cqfd\\

\textit{Proof of Theorem \ref{MITThmConv2}}.  (a) Direct calculus gives 
\begin{equation}\label{derivee 1 nulle}
 \P f'_{\bo} - \P^\perp g'_{\bo} =0
\end{equation}
and 
\begin{equation}\label{derivee 2 negative}
 \P f''_{\bo} - \P^\perp g''_{\bo} = -\P^\perp \left(h_{\bo}' h_{\bo}'^\top \right) = -\Sigma_1.
\end{equation}
Observe that the above matrix $\Sigma_1$ is symmetric and positive.\\
For any $\theta\in\Theta$, we have  $M_n'(\theta) = \widehat{\P} f'_{\theta} - \widehat{\P}^\perp g'_{\theta}.$ Note that  $$f'_{\bo}(x,y) = g'_{\bo}(x,y)=
\left(1,\xi_1(x)\zeta_1(y),\ldots,\xi_d(x)\zeta_d(y)\right)^\top.$$
We will state the asymptotic normality of $\sqrt{n} M_n'(\boldsymbol{0})$ using the multivariate Delta method. 
So consider the random column vector in $\R^{1+3d}$
$$V(X,Y):=\left(1, \xi_1(X),\ldots,\xi_d(X),\zeta_1(Y),\ldots,\zeta_d(Y),\xi_1(X)\zeta_1(Y),\ldots,\xi_d(X)\zeta_d(Y)\right)^\top.$$
Denote by 
$$\mu := \E(V(X,Y)) = \left( 1, \P_1 \xi_1,\ldots,\P_1 \xi_d,\P_2 \zeta_1,\ldots,\P_2 \zeta_d, \P_1 \xi_1 \P_2 \zeta_1,\ldots,\P_1 \xi_d\P_2 \zeta_d\right)^
\top$$ 
which is a column vector in $\R^{1+3d}$. Then we have, by multivariate central limit theorem, the convergence in distribution   
$$\sqrt{n}\left(\frac{1}{n}\sum_{i=1}^n V(X_i,Y_i) - \mu\right) \to \mathcal{N}_{1+3d}\left(0,\Sigma\right),$$
with $\Sigma = \E\left((V(X,Y)-\mu)(V(X,Y)-\mu)^\top\right)$,
from which we obtain, by multivariate Delta method,
\begin{equation} \label{MITEqnVar2}
 \sqrt{n}\left(M_n'(\boldsymbol{0}) - \psi(\mu)\right) \to \mathcal{N}_{1+d}\left(\bo, \Sigma_2 := \psi'(\mu) \Sigma \psi'(\mu)^\top\right),
\end{equation}
where $\psi(\cdot)$ is the function defined on $\R^{1+3d}$ into $\R^{1+d}$ by $$\psi(x_0,x_1,\ldots,x_d,y_1,\ldots,y_d,z_1,\ldots,z_d) = (0,x_1y_1-
z_1,\ldots,x_dy_d-z_d)^\top$$
which is of class $\mathcal{C}^1.$ Note that $\psi(\mu)=\bo$, the first component of $M_n'(\bo)$ is equal to zero for all $n$
and that the first column and row of the limiting 
covariance matrix $\Sigma_2$ are equal both to $\bo$. Whence we have the convergence in distribution
\begin{equation}\label{TCL 1}
\sqrt{n} M_n'(\boldsymbol{0}) \to \mathcal{N}_{1+d}\left(\bo, \Sigma_2\right).
\end{equation}
By Taylor expansion of  $U_n(\widehat{\theta}_{\varphi_1})$ in $\widehat{\theta}_{\varphi_1}$ around $\theta_T=\bo$, using condition (A.8)
and the convergence in probability of  $\widehat{\theta}_{\varphi_1}$ to $\theta_T=\bo$, we obtain
\begin{equation}\label{dev Tay 1}
\bo = M_n'(\widehat{\theta}_{\varphi_1}) = M_n'(\bo)+M_n''(\bo) \widehat{\theta}_{\varphi_1} +o_\P(1) \widehat{\theta}_{\varphi_1}.
\end{equation}
On the other hand, by (A.9), we can write 
$$ M_n''(\bo) = \P f''_{\bo} - \P^\perp g''_{\bo} +o_\P(1) = - \Sigma_1 + o_\P(1).$$
Combining the last two displays, leads to 
\begin{equation} \label{vitesse 1}
  M_n'(\bo) = \left(\Sigma_1+o_\P(1)\right)  \widehat{\theta}_{\varphi_1}.
\end{equation}
We have, from (\ref{TCL 1}),  that $\sqrt{n} M_n'(\boldsymbol{0})=O_\P(1)$, which by (\ref{vitesse 1}) implies that 
$\sqrt{n}\widehat{\theta}_{\varphi_1} =O_\P(1).$
Combining this last result with the relation (\ref{dev Tay 1}), we obtain 
\begin{equation}\label{dev 2 de theta} 
\sqrt{n}\widehat{\theta}_{\varphi_1} = \Sigma_1^{-1} \sqrt{n} M_n'(\bo)+o_\P(1).
\end{equation}
Use this last relation and (\ref{TCL 1}) to conclude the proof of part (a).\\

(b) By Taylor expansion of $\widehat{I}_{\varphi_1} =  M_n(\widehat{\theta}_{\varphi_1})$, in $\widehat{\theta}_{\varphi_1}$
around $\theta_T=\bo$, using the fact that $M_n(\bo)=0$ and some of the above statements, we obtain 
\begin{eqnarray}
 \widehat{I}_{\varphi_1} & := & M_n(\widehat{\theta}_{\varphi_1}) \nonumber\\
   & = &  M_n'(\bo) \widehat{\theta}_{\varphi_1} - \frac{1}{2}  \widehat{\theta}_{\varphi_1}^\top \Sigma_1 
   \widehat{\theta}_{\varphi_1} +o_\P(1/n)\nonumber
\end{eqnarray}
which by (\ref{dev 2 de theta}) leads to 
\begin{eqnarray}\label{loi la stat km}
2 n \, \widehat{I}_{\varphi_1} & = & \left(\sqrt{n}M_n'(\bo)\right)^\top  \Sigma_1^{-1} \sqrt{n}M_n'(\bo) +o_\P(1)\\
     & = & \left(\sqrt{n} \Sigma_1^{-1/2} M_n'(\bo)\right)^\top  \Sigma_1^{-1/2} \sqrt{n}M_n'(\bo) +o_\P(1).
\end{eqnarray}
This proves the convergence in distribution of $2 n \, \widehat{I}_{\varphi_1}$ to the random variable $Z^\top Z$,
where $Z$ is a centered multivariate normal random variable with covariance matrix $C = \Sigma_1^{-1/2} \Sigma_2 \Sigma_1^{-1/2}.$
\cqfd\\

\textit{Proof of Theorem \ref{LDP}}.
First, under assumption (A.10),
\cite{Eischelsbacher_Scmock2002_AIHP} yields the following large deviations principle for the sequence $(\widehat{\P})_n$ of empirical measures :  we have 
for all measurable subset $B$ of $\mathcal{M}_\mathcal{G}$, 
\begin{equation} \label{BEITeqnLDPech1}
 \liminf_{n \rightarrow \infty} \frac{1}{n} \log \P^{\perp} \left(\widehat{\P} \in B\right) 
\,  \geq \,  - \inf_{Q\in \text{ Int}_{\tau_\mathcal{G}}(B)}\mathbb{K}(Q,\P^{\perp}), 
  \end{equation}
  and 
  \begin{equation} \label{BEITeqnLDPech2}
 \limsup_{n \rightarrow \infty} \frac{1}{n} \log \P^{\perp} \left(\widehat{\P} \in B\right) 
 \,  \leq \,  - \inf_{Q \in \text{ Cl}_{\tau_\mathcal{G}} (B)} \mathbb{K}(Q,\P^{\perp}),
\end{equation}
where $\mbox{Int}_{\tau_\mathcal{G}} (B)$ and $\mbox{Cl}_{\tau_\mathcal{G}} (B)$ 
denote, respectively, the interior and closure of $B$, with respect to the $\tau_\mathcal{G}$-topology.
Since $Q\in\mathcal{M}_\mathcal{G}\mapsto 
\Ddp(Q,Q^{\perp})$ is continuous, we obtain by contraction principle from (\ref{BEITeqnLDPech1}) and (\ref{BEITeqnLDPech2}), for all $d>0$, 
\begin{equation}  \label{BEITeqnLDPDinf}
 \liminf_{n \rightarrow \infty} \frac{1}{n} \log \P^{\perp} \left( \widehat{I_{\varphi}} >d \right)  \, \geq \,  - \inf\left\{\mathbb{K}(Q,\P^{\perp}); \, 
 Q\in\mathcal{M}_\mathcal{G} \text{ and } \Ddp(Q,Q^{\perp}) >d\right\}
 \end{equation}
 and 
\begin{equation}   \label{BEITeqnLDPsup}
 \limsup_{n \rightarrow \infty} \frac{1}{n} \log \P^{\perp} \left(\widehat{I_{\varphi}} >d\right)  \, \leq \, - \inf\left\{\mathbb{K}(Q,\P^{\perp}); \, 
 Q\in\mathcal{M}_\mathcal{G} \text{ and } \Ddp(Q,Q^{\perp}) \geq d\right\}.
\end{equation}
We now prove that the function $e_{\varphi}(\cdot) : d \in \mathbb{R}_+^* \mapsto \inf \{\mathbb{K} (Q, \P^\perp); \, Q\in\mathcal{M}_\mathcal{G} \text{ and } 
\Ddp(Q,Q^{\perp}) \geq d\} \in [0,+\infty]$ is right-continuous so that infima in 
(\ref{BEITeqnLDPDinf}) and (\ref{BEITeqnLDPsup}) are equal, yielding (\ref{BEITEqnLDP}). 
So, let $d>0$ be any positive real number, and show that $e_\varphi(\cdot)$ is right-continuous at $d$. 
If no $Q \in \Omega_d$ exists such that $\K(Q,\P^\perp) < + \infty$,
obviously, since for any $d' \in \R_+^*$ such that $d \leq d'$, we have $\Omega_{d'} \subseteq \Omega_d$, 
then both $e_\varphi(d)$ and $e_\varphi(d')$ equal $\infty$, which implies that $e_\varphi$ is 
right-continuous at $d$ in this case. Now, assume that some $Q \in \Omega_d$ exists such that $\K(Q,\P^\perp) < \infty$. 
Two cases can be handled 
separately.
First, assume that the infimum~(\ref{BEITeqnFTLDP}) is achieved for some $Q$ such that $\mathcal{D}_\varphi(Q, Q^\perp)=:d'>d$. Then, for all $d''$ satisfying $d\leq 
d''\leq d'$, the equality $e_\varphi(d'') = e_{\varphi}(d)$ holds, yielding the right-continuity of $e_\varphi$ at $d$.
Second, assume that the infimum~(\ref{BEITeqnFTLDP}) is achieved for $Q$ such that $\mathcal{D}_\varphi(Q, Q^\perp)=d$.
Let us prove that there exists a sequence $(Q_n)_n$ of elements of $\{Q : \mathcal{D}_\varphi(Q, Q^\perp) >d \}$ such that $\K(Q_n, \P^\perp) \stackrel{n 
\rightarrow \infty}{\longrightarrow} \mathbb{K}(Q,\P^\perp)$  yielding right-continuity of $e_\varphi(\cdot)$ at  $d$.
We build such a sequence $(Q_n)_n$ such that $Q_n$ has the same marginal distributions as  
$Q$, i.e., $Q_{n,1} =Q_1$ and $Q_{n,2} = Q_2$. We have then $Q_n^{\perp} = Q^{\perp}$.
Let 
\begin{equation*}
\overline{\theta} := \arg\sup_{\theta \in \Theta} \left\{ \int \varphi'( h_{\theta}) \, dQ - \int \varphi^* ( \varphi' ( h_{\theta} )) \, dQ^{\perp} \right\},
\end{equation*}
so that
\begin{equation} \label{MITEqnPseudoDiv}
 \int \varphi'(h_{\overline{\theta}}) \, dQ - \int \varphi^* \left(\varphi'(h_{\overline{\theta}})\right)\,  dQ^{\perp} = \Ddp(Q,Q^{\perp}) = d>0.
\end{equation}
Denote $\widetilde{Q}$ the image distribution on the Borel  $\sigma$-field $\left(\R,\mathcal{B}(\R)\right)$ of $Q$ by the function $\varphi' ( h_{\overline{\theta}} )$. 
Let us prove by contradiction  that $\widetilde{Q}$ can not be Dirac measure, by making use of either (A.12.a) or (A.12.b). If $\widetilde{Q}$ was a Dirac measure,
 necessarily 
$\varphi'(h_{\overline{\theta}})$ would be $Q$-a.s. constant, i.e., $h_{\overline{\theta}}$ would be $Q$-a.s. constant 
\begin{equation} \label{MITEqnProofLDP1_1}
 h_{\overline{\theta}}(\cdot,\cdot) = c, \quad Q\textrm{-a.s.}
\end{equation}
Now, if (A.12.a) holds, we can consider the set of all finite-discrete distributions with the same finite 
support $\Om$, instead of the set $\mathcal{M}_\mathcal{G}$. 
Hence,  $Q$ and $Q^\perp$ have same support, so that (\ref{MITEqnProofLDP1_1}) implies that 
\begin{equation} \label{MITEqnProofLDP1_2}
 h_{\overline{\theta}}(\cdot,\cdot) = c, \quad Q^\top\textrm{-a.s.}
\end{equation}
Combining (\ref{MITEqnProofLDP1_1}), (\ref{MITEqnProofLDP1_2}) and (\ref{MITEqnPseudoDiv}), 
we obtain
\begin{equation}\label{relation between c and d}
\varphi(c)+\varphi'(c)(1-c)=d>0.
\end{equation}
On the other hand, by convexity of $\varphi(\cdot)$ and the fact that $\varphi(1)=0$, we get
$$0=\varphi(1)\geq \varphi(c)+\varphi'(c)(1-c)=d,$$   
which contradicts the fact that $d>0$.
Alternatively, assume that (A.12.b) holds. 
Note that, under this assumption in connection with (A.2), we can see that the value $\theta_0$ (of the parameter
corresponding to independence) is necessarly $\theta_0 :=(\alpha_0,\beta_0^\top)^\top =(0,0,\ldots, 0)^\top$.
We can see also, by contradiction as above, that $\overline{\theta}$ can not be of the form $(\overline{\alpha},0,\ldots,0)^\top$ 
with $\overline{\alpha}\neq 0$. Hence, it can be written as 
\begin{equation}\label{theta bar et beta bar}
\overline{\theta}=(\overline{\alpha},\overline{\beta}^\top)^\top \quad \text{with} \quad
\overline{\beta}\neq (0,\ldots,0)^\top.
\end{equation} 
Now, by (\ref{MITEqnProofLDP1_1}), using the fact that $h_\theta(\cdot,\cdot)$ is of the form $\exp\left(\alpha+m_\beta(\cdot,\cdot)\right)$,
we get that  
\begin{equation} \label{MITEqnProofLDP1_3}
m_{\overline{\beta}}(\cdot,\cdot) = cte, \quad Q\textrm{-a.s.}
\end{equation}
Note that the support of $Q^\perp$ is included in that of $\P^\perp$ (if not, $Q$ would not be a.c.w.r.t. $\P^\perp$ and 
$\K(Q, \P^\perp)$ would not be finite). Hence, (\ref{MITEqnProofLDP1_3}) 
implies that $\P^\perp(m_{\overline{\beta}}(X,Y) =cte) \neq 0$, which in turn implies that $\overline{\beta} =(0,\ldots,0)^\top$ by  assumption (A.12.b). 
This contradicts (\ref{theta bar et beta bar}).
We have proven then that $\widetilde{Q}$ is not a Dirac measure. 
So, there exist $A$, $B$ two measurable subsets of $(a_{\varphi^*}, b_{\varphi^*}) =\textrm{Im}
(\varphi')$ 
such that $\widetilde{Q}(A) >  0$, $\widetilde{Q}
(B) > 0$ and $a:=\inf A > b:=\sup(B)$.
Denoting $g := dQ/ dQ^{\perp}$ the density of $Q$ with respect to the product of its marginal distributions, set
\begin{eqnarray*}
 g_n & := & \left( 1+ \frac{c_1}{n}\right) g \id_{\{\varphi' (h_{\overline{\theta}}) \in A\}} + \left( 1- \frac{c_2}{n}\right) g \id_{\{\varphi'(h_{\overline{\theta}})\in 
 B\}} + g \id_{\{\varphi' (h_{\overline{\theta}}) \in \overline{A \cup B}\}} \\
  & = & g+ \frac{c_1}{n}g \id_{\{\varphi' (h_{\overline{\theta}}) \in A\}} - \frac{c_2}{n} g \id_{\{\varphi'(h_{\overline{\theta}} )\in B\}},
\end{eqnarray*}
where $c_1 := \widetilde{Q}(B)$ and $c_2 := \widetilde{Q}(A)$. 
Note that $g_n$ is nonnegative for $n$ sufficiently large, and that $\int_{\Om}g_n(x,y)\, dQ^\perp(x,y)=1$.
Then, let $Q_n$ be the probability distribution on $\Om$ such that $Q_{n,1} = Q_1$, 
$Q_{n,2} = Q_2$ and $\der Q_n/ \der Q^{\perp} = g_n$. 
We have 
\begin{eqnarray*}
 \int_{\Om} \varphi'( h_{\overline{\theta}}) \, dQ_n & = & \int_{\Om} \varphi'( h_{\overline{\theta}}) g_n\, dQ^{\perp} \\
  & = & \int \varphi'( h_{\overline{\theta}})
  \,dQ + \frac{c_1}{n} \E_{\widetilde{Q}} (\mbox{Id}.\id_{A}) - \frac{c_2}{n} \E_{\widetilde{Q}} (\mbox{Id}.\id_{B}) \\
  & \geq &  \int \varphi' (h_{\overline{\theta}})\, dQ +  \frac{c_1}{n} a \widetilde{Q}(A) - \frac{c_2}{n}b \widetilde{Q}(B) \\
  & = &  \int \varphi' (h_{\overline{\theta}})\, dQ + \frac{c_1 c_2}{n} (a-b) \\
  & > & \int \varphi' (h_{\overline{\theta}})\, dQ,
\end{eqnarray*}
where $\mbox{Id}(x) := x$, for all $x \in (a_{\varphi^*}, b_{\varphi^*})$.
Then, 
\begin{eqnarray*}
 \Ddp(Q_n,Q_n^{\perp}) & = & \int \varphi' (h_{\overline{\theta}})\, dQ_n  -  \int \varphi^* ( \varphi' ( h_{\overline{\theta}}))\, dQ^{\perp} \\
  & > & \int \varphi' ( h_{\overline{\theta}})\, dQ -  \int \varphi^* (\varphi' (h_{\overline{\theta}}))\, dQ^{\perp} \\
   & = & d.
\end{eqnarray*}
Finally, the convergence of $\mathbb{K}(Q_n,\P^\perp)$ to $\mathbb{K}(Q,\P^\perp)$ can be proved using the decompositions
 $$\mathbb{K}(Q_n,\P^\perp) = \mathbb{K}(Q_n,Q^\perp)+\mathbb{K}(Q^\top,\P^\perp), \quad \mathbb{K}(Q,\P^\perp) = \mathbb{K}(Q,Q^\perp)+
 \mathbb{K}(Q^\top,\P^\perp),$$
and Lebesgue's dominated convergence theorem.
\cqfd

\bibliographystyle{natbib}

\end{document}